\documentclass{amsart}%
\usepackage{amsfonts}
\usepackage{graphicx}
\usepackage{amscd}%
\usepackage{amsmath}%
\usepackage{amssymb}
\newtheorem{theorem}{Theorem}
\theoremstyle{plain}

\newtheorem{case}{Case}

\newtheorem{corollary}{Corollary}

\newtheorem{proposition}{Proposition}
\newtheorem{remark}{Remark}

\numberwithin{equation}{section}
\numberwithin{theorem}{section}
\numberwithin{algorithm}{section}
\numberwithin{axiom}{section}
\numberwithin{case}{section}
\numberwithin{claim}{section}
\numberwithin{conclusion}{section}
\numberwithin{condition}{section}
\numberwithin{conjecture}{section}
\numberwithin{corollary}{section}
\numberwithin{criterion}{section}
\numberwithin{definition}{section}
\numberwithin{example}{section}
\numberwithin{exercise}{section}
\numberwithin{lemma}{section}
\numberwithin{notation}{section}
\numberwithin{problem}{section}
\numberwithin{proposition}{section}
\numberwithin{remark}{section}
\numberwithin{solution}{section}

\begin{document}
\title{Sharp integral inequalities for harmonic functions}
\author{Fengbo Hang}
\address{Department of Mathematics\\
Princeton University\\
Fine Hall, Washington Road\\
Princeton, NJ 08544}
\email{fhang@math.princeton.edu}
\author{Xiaodong Wang}
\address{Department of Mathematics\\
Michigan State University\\
East Lansing, MI 48824}
\email{xwang@math.msu.edu}
\author{Xiaodong Yan}
\address{Department of Mathematics\\
Michigan State University\\
East Lansing, MI 48824}
\email{xiayan@math.msu.edu}

\begin{abstract}
Motivated by Carleman's proof of isoperimetric inequality in the plane, we
study some sharp integral inequalities for harmonic functions on the upper
halfspace. We also derive the regularity for nonnegative solutions of the
associated integral system and some Liouville type theorems.

\end{abstract}
\maketitle

\section{Introduction\label{sec1}}

The classical isoperimetric inequality in the plane states that for any
bounded domain with area $A$ and boundary length $L$ we have%
\begin{equation}
4\pi A\leq L^{2}\label{eq1.1}%
\end{equation}
and equality holds if and only if the domain is a disk. Inequality
(\ref{eq1.1}) remains true for bounded domains in a simply connected surface
with nonpositive curvature. Among the proofs of this fact the one due to
Carleman \cite{C} is particularly interesting. Indeed, let $\left(
M^{2},g\right)  $ be any simply connected compact surface with boundary and
nonpositive curvature, it follows from Riemann mapping theorem that $\left(
M^{2},g\right)  $ is isometric to $\left(  \overline{B}_{1}^{2},e^{2w}%
g_{\mathbb{R}^{2}}\right)  $, here $B_{1}^{2}$ is the two dimensional open
unit disk and $g_{\mathbb{R}^{2}}$ is the Euclidean metric on $\mathbb{R}^{2}%
$. The nonpositivity of curvature implies that $w$ is a subharmonic function.
Let $u$ be the harmonic function on $B_{1}$ with the same boundary value as
$w$, then $w\leq u$. In \cite{C} it was proved that for any smooth harmonic
function on $\overline{B}_{1}^{2}$ we have%
\begin{equation}
\int_{B_{1}}e^{2u}dx\leq\frac{1}{4\pi}\left(  \int_{S^{1}}e^{u}d\theta\right)
^{2}\label{eq1.2}%
\end{equation}
and equality holds if and only if $u\left(  x\right)  =c$ or $-2\log\left\vert
x-x_{0}\right\vert +c$ for some $x_{0}\in\mathbb{R}^{2}\backslash\overline
{B}_{1}$ and constant $c$. We may ask for natural generalizations to higher
dimensions. Without an analog of the Riemann mapping theorem, we may start
with a metric $g=\rho^{\frac{4}{n-2}}g_{\mathbb{R}^{n}}$ on $\overline{B}%
_{1}^{n}$ with nonpositive scalar curvature, here $n\geq3$, $B_{1}^{n}$ is the
open unit ball in $\mathbb{R}^{n}$ and $g_{\mathbb{R}^{n}}$ is the Euclidean
metric on $\mathbb{R}^{n}$. It follows that $\rho$ is a subharmonic function.
Under the metric $g$ the volume of $\overline{B}_{1}$ is equal to $\int
_{B_{1}}\rho^{\frac{2n}{n-2}}dx$ and the area of $\partial B_{1}$ is equal to
$\int_{\partial B_{1}}\rho^{\frac{2\left(  n-1\right)  }{n-2}}dS$. We would
like to know whether the inequality%
\[
\int_{B_{1}}\rho^{\frac{2n}{n-2}}dx\leq n^{-\frac{n}{n-1}}\omega_{n}%
^{-\frac{1}{n-1}}\left(  \int_{\partial B_{1}}\rho^{\frac{2\left(  n-1\right)
}{n-2}}dS\right)  ^{\frac{n}{n-1}}%
\]
is still true. Here $\omega_{n}$ is the Euclidean volume of the unit ball in
$\mathbb{R}^{n}$. Since $\rho$ is bounded from above by the harmonic function
with the same boundary value, we only need to know whether the inequality%
\begin{equation}
\left\vert u\right\vert _{L^{\frac{2n}{n-2}}\left(  B_{1}\right)  }\leq
n^{-\frac{n-2}{2\left(  n-1\right)  }}\omega_{n}^{-\frac{n-2}{2n\left(
n-1\right)  }}\left\vert u\right\vert _{L^{\frac{2\left(  n-1\right)  }{n-2}%
}\left(  \partial B_{1}\right)  }\label{eq1.3}%
\end{equation}
is true for every smooth harmonic function $u$ on $\overline{B}_{1}^{n}$. The
answer to this question is affirmative and the inequality may be proved by
subcritical approximation (see \cite{HWY}). However, for future purpose it
seems helpful to transfer this problem to upper halfspace and derive some
Liouville type results. Indeed, assume $u$ is a positive harmonic function on
$\overline{B}_{1}$, let $e_{n}=\left(  0,\cdots,0,1\right)  $ and $\phi$ be
the Mobius transformation given by%
\[
\phi\left(  x\right)  =\frac{x+\frac{e_{n}}{2}}{\left\vert x+\frac{e_{n}}%
{2}\right\vert ^{2}}-e_{n}.
\]
Then $\phi\left(  \mathbb{R}_{+}^{n}\right)  =B_{1}$ and%
\[
\phi^{\ast}g_{\mathbb{R}^{n}}=\frac{1}{\left\vert x+\frac{e_{n}}{2}\right\vert
^{4}}\sum_{i=1}^{n}dx_{i}\otimes dx_{i}.
\]
Here $\mathbb{R}_{+}^{n}=\left\{  x\in\mathbb{R}^{n}:x_{n}>0\right\}  $. Let%
\[
v\left(  x\right)  =\frac{1}{\left\vert x+\frac{e_{n}}{2}\right\vert ^{n-2}%
}u\left(  \frac{x+\frac{e_{n}}{2}}{\left\vert x+\frac{e_{n}}{2}\right\vert
^{2}}-e_{n}\right)  ,
\]
then $\phi^{\ast}\left(  u^{\frac{4}{n-2}}g_{\mathbb{R}^{n}}\right)
=v^{\frac{4}{n-2}}g_{\mathbb{R}^{n}}$. The inequality (\ref{eq1.3}) becomes%
\begin{equation}
\left\vert v\right\vert _{L^{\frac{2n}{n-2}}\left(  \mathbb{R}_{+}^{n}\right)
}\leq n^{-\frac{n-2}{2\left(  n-1\right)  }}\omega_{n}^{-\frac{n-2}{2n\left(
n-1\right)  }}\left\vert v\right\vert _{L^{\frac{2\left(  n-1\right)  }{n-2}%
}\left(  \mathbb{R}^{n-1}\right)  }.\label{eq1.4}%
\end{equation}
Note that since $v$ is the Poission integral of $\left.  v\right\vert
_{\mathbb{R}^{n-1}}$, inequality (\ref{eq1.4}) follows from Theorem
\ref{thm1.1} below. To state the results, let us fix some notations. For
convenience, we use $x,y,\cdots$ for points in $\mathbb{R}^{n}$ and $\xi
,\zeta,\cdots$ for points in $\mathbb{R}^{n-1}=\left\{  x\in\mathbb{R}%
^{n}:x_{n}=0\right\}  $. For $x\in\mathbb{R}^{n}$, we let $x^{\prime}=\left(
x_{1},\cdots,x_{n-1}\right)  $, $x=\left(  x^{\prime},x_{n}\right)  $. The
Poission kernel for the upper half space is given by (see \cite[p61]{S})%
\[
P\left(  x,\xi\right)  =\frac{2}{n\omega_{n}}\frac{x_{n}}{\left(  \left\vert
x^{\prime}-\xi\right\vert ^{2}+x_{n}^{2}\right)  ^{n/2}}\text{\quad for }%
x\in\mathbb{R}_{+}^{n}\text{, }\xi\in\mathbb{R}^{n-1}\text{.}%
\]
Given a function $f$ defined on $\mathbb{R}^{n-1}$, let%
\[
\left(  Pf\right)  \left(  x\right)  =\int_{\mathbb{R}^{n-1}}P\left(
x,\xi\right)  f\left(  \xi\right)  d\xi\text{\quad for }x\in\mathbb{R}_{+}%
^{n}.
\]
We have the following sharp inequality for $P$ (see Theorem \ref{thm4.1}):

\begin{theorem}
\label{thm1.1}Assume $n\geq3$, then for any $f\in L^{\frac{2\left(
n-1\right)  }{n-2}}\left(  \mathbb{R}^{n-1}\right)  $,%
\begin{equation}
\left\vert Pf\right\vert _{L^{\frac{2n}{n-2}}\left(  \mathbb{R}_{+}%
^{n}\right)  }\leq n^{-\frac{n-2}{2\left(  n-1\right)  }}\omega_{n}%
^{-\frac{n-2}{2n\left(  n-1\right)  }}\left\vert f\right\vert _{L^{\frac
{2\left(  n-1\right)  }{n-2}}\left(  \mathbb{R}^{n-1}\right)  }.\label{eq1.5}%
\end{equation}
Moreover, equality holds if and only if $f\left(  \xi\right)  =\frac
{c}{\left(  \lambda^{2}+\left\vert \xi-\xi_{0}\right\vert ^{2}\right)
^{\frac{n-2}{2}}}$ for some constant $c$, positive constant $\lambda$ and
$\xi_{0}\in\mathbb{R}^{n-1}$.
\end{theorem}

If we look at the variational problem%
\begin{equation}
c_{n}=\sup\left\{  \left\vert Pf\right\vert _{L^{\frac{2n}{n-2}}\left(
\mathbb{R}_{+}^{n}\right)  }:f\in L^{\frac{2\left(  n-1\right)  }{n-2}}\left(
\mathbb{R}^{n-1}\right)  ,\left\vert f\right\vert _{L^{\frac{2\left(
n-1\right)  }{n-2}}\left(  \mathbb{R}^{n-1}\right)  }=1\right\}
.\label{eq1.6}%
\end{equation}
Then any nonnegative critical function $f$, after scaling must satisfy%
\begin{equation}
f\left(  \xi\right)  ^{\frac{n}{n-2}}=\int_{\mathbb{R}_{+}^{n}}P\left(
x,\xi\right)  \left(  Pf\right)  \left(  x\right)  ^{\frac{n+2}{n-2}%
}dx.\label{eq1.7}%
\end{equation}
We have the following Liouville type theorem (see Proposition \ref{prop6.1})
which is in the same spirit as those in \cite{GNN,CGS,CLO2,L}.

\begin{theorem}
\label{thm1.2}Assume $n\geq3$, $f\in L_{loc}^{\frac{2\left(  n-1\right)
}{n-2}}\left(  \mathbb{R}^{n-1}\right)  $ is nonnegative, not identically zero
and it satisfies (\ref{eq1.7}), then for some $\lambda>0$ and $\xi_{0}%
\in\mathbb{R}^{n-1}$,%
\[
f\left(  \xi\right)  =c\left(  n\right)  \left(  \frac{\lambda}{\lambda
^{2}+\left\vert \xi-\xi_{0}\right\vert ^{2}}\right)  ^{\frac{n-2}{2}}.
\]

\end{theorem}

We note that the condition $f\in L_{loc}^{\frac{2\left(  n-1\right)  }{n-2}%
}\left(  \mathbb{R}^{n-1}\right)  $ can not be dropped since $c\left(
n\right)  \left\vert \xi\right\vert ^{-\frac{n-2}{2}}$ is a singular solution
for (\ref{eq1.7}). During the process of identifying maximizing functions in
Theorem \ref{thm1.1} and the critical functions in Theorem \ref{thm1.2}, we
establish the following interesting fact (see Proposition \ref{prop4.1}):

\begin{proposition}
\label{prop1.1}Let $n\geq2$, $u$ be a function on $\mathbb{R}^{n}$ which is
radial with respect to the origin, $0<u\left(  x\right)  <\infty$ for $x\neq0
$, $e_{1}=\left(  1,0,\cdots,0\right)  $, $\alpha\in\mathbb{R}$, $\alpha\neq
0$. If $v\left(  x\right)  =\left\vert x\right\vert ^{\alpha}u\left(  \frac
{x}{\left\vert x\right\vert ^{2}}-e_{1}\right)  $ is radial with respect to
some point, then either $u\left(  x\right)  =\left(  c_{1}\left\vert
x\right\vert ^{2}+c_{2}\right)  ^{\alpha/2}$ for some $c_{1}\geq0$, $c_{2}>0 $
or%
\[
u\left(  x\right)  =\left\{
\begin{tabular}
[c]{l}%
$c_{1}\left\vert x\right\vert ^{\alpha},$ if $x\neq0,$\\
$c_{2},$ if $x=0,$%
\end{tabular}
\right.
\]
for some $c_{1}>0$ and $c_{2}$, an arbitrary number.
\end{proposition}

There are similar statements for the cases $\alpha=0$ or $n=1$ (see Remark
\ref{rmk4.1} and Proposition \ref{prop4.2}). The crucial point of Proposition
\ref{prop1.1} is that we do not need any regularity assumption on the function
$u$. This is very convenient when the regularity of extremal functions are
hard to get apriorly. The radial symmetry property of function may come from
symmetrization arguments or the method of moving planes etc. For example,
Proposition \ref{prop1.1} gives another way to determine the maximizing
functions for those cases of Hardy-Littlewood-Sobolev inequalities studied in
\cite[section III]{Li2}. The formulation of Proposition \ref{prop1.1} is
motivated from previous works in \cite{CL,O}, \cite[section 3]{CLO2} and
\cite[section 6]{CLO3}. It is worth pointing out that Proposition
\ref{prop1.1} is the fact for method of moving planes which corresponds to the
fact \cite[lemma 2.5]{LZ} or \cite[lemma 5.8]{L} for the method of moving
spheres, a variant of the method of moving planes.

According to Proposition \ref{propinterpolation} below, for $n\geq2$ and
$1<p<\infty$ the operator%
\[
P:L^{p}\left(  \mathbb{R}^{n-1}\right)  \rightarrow L^{\frac{np}{n-1}}\left(
\mathbb{R}_{+}^{n}\right)  :f\mapsto Pf
\]
is always a bounded linear map. From the analytical point view it is
interesting to consider the variational problem%
\begin{equation}
c_{n,p}=\sup\left\{  \left\vert Pf\right\vert _{L^{\frac{np}{n-1}}\left(
\mathbb{R}_{+}^{n}\right)  }:f\in L^{p}\left(  \mathbb{R}^{n-1}\right)
,\left\vert f\right\vert _{L^{p}\left(  \mathbb{R}^{n-1}\right)  }=1\right\}
\label{eq1.8}%
\end{equation}
for all such $p$'s. Fix $1<p<\infty$, for a function $f$ defined on
$\mathbb{R}^{n-1}$, $\lambda>0$ and $\zeta\in\mathbb{R}^{n-1}$, we write%
\[
f^{\lambda,\zeta}\left(  \xi\right)  =\lambda^{-\frac{n-1}{p}}f\left(
\frac{\xi-\zeta}{\lambda}\right)  \text{\quad for }\xi\in\mathbb{R}^{n-1}.
\]
Then we have (see Theorem \ref{thm3.1} and Theorem \ref{thm4.1}):

\begin{theorem}
\label{thm1.3}Given $n\geq2$ and $1<p<\infty$.

\begin{itemize}
\item Let $f_{i}$ be a maximizing sequence of functions for (\ref{eq1.8}),
then after passing to a subsequence there exists $\lambda_{i}>0$ and
$\zeta_{i}\in\mathbb{R}^{n-1}$ such that $f_{i}^{\lambda_{i},\zeta_{i}%
}\rightarrow f$ in $L^{p}\left(  \mathbb{R}^{n-1}\right)  $. In particular,
there exists at least one maximizing function.

\item After multiplying by a nonzero constant, every maximizer $f$ of
(\ref{eq1.8}) is nonnegative, radial symmetric with respect to some point,
strictly decreasing in the radial direction and it satisfies%
\begin{equation}
f\left(  \xi\right)  ^{p-1}=\int_{\mathbb{R}_{+}^{n}}P\left(  x,\xi\right)
\left(  Pf\right)  \left(  x\right)  ^{\frac{np}{n-1}-1}dx.\label{eq1.9}%
\end{equation}

\item If $n\geq3$ and $p=\frac{2\left(  n-1\right)  }{n-2}$, then any
maximizer of (\ref{eq1.8}) must be of the form%
\[
f\left(  \xi\right)  =\pm c\left(  n\right)  \left(  \frac{\lambda}%
{\lambda^{2}+\left\vert \xi-\xi_{0}\right\vert ^{2}}\right)  ^{\frac{n-2}{2}}%
\]
for some $\lambda>0$, $\xi_{0}\in\mathbb{R}^{n-1}$. In particular
$c_{n,\frac{2\left(  n-1\right)  }{n-2}}=n^{-\frac{n-2}{2\left(  n-1\right)
}}\omega_{n}^{-\frac{n-2}{2n\left(  n-1\right)  }}$.

\item If $n\geq3$ and $p=\frac{2\left(  n-1\right)  }{n}$, then any maximizer
of (\ref{eq1.8}) must be of the form%
\[
f\left(  \xi\right)  =\pm c\left(  n\right)  \left(  \frac{\lambda}%
{\lambda^{2}+\left\vert \xi-\xi_{0}\right\vert ^{2}}\right)  ^{n/2}%
\]
for some $\lambda>0$, $\xi_{0}\in\mathbb{R}^{n-1}$. In particular
\[
c_{n,\frac{2\left(  n-1\right)  }{n}}=\frac{1}{\sqrt{2\left(  n-2\right)
}\sqrt[4]{\pi}}\left(  \frac{\left(  n-2\right)  !}{\Gamma\left(  \frac
{n-1}{2}\right)  }\right)  ^{\frac{1}{2\left(  n-1\right)  }}.
\]

\end{itemize}
\end{theorem}

It is interesting that the problem considered here demonstrates very similar
structures to the sharp Hardy-Littlewood-Sobolev inequalities studied in
\cite{Li2}. Besides above properties of maximizing functions, we know they are
smooth. This is a nontrivial fact since it does not follow from the usual
bootstrap method. Indeed, we know all the nonnegative critical functions of
(\ref{eq1.8}) are smooth and radial symmetric with respect to some point (see
Theorem \ref{thm5.1} and Theorem \ref{thm6.1}). More precisely we have

\begin{theorem}
\label{thm1.4}Given $n\geq2$ and $1<p<\infty$. If $f\in L^{p}\left(
\mathbb{R}^{n-1}\right)  $ is nonnegative, not identically zero and it
satisfies (\ref{eq1.9}), then $f\in C^{\infty}\left(  \mathbb{R}^{n-1}\right)
$, moreover it is radial symmetric with respect to some point and strictly
decreasing along the radial direction.
\end{theorem}

In Section \ref{sec2} below, we will collect some basic estimates for Poission
integrals and show the operator $P$ is bounded in suitable Lebesgue spaces and
Lorentz spaces. In Section \ref{sec3}, we apply the general frame of
concentration compactness principle (\cite{Lion}) to show that every
maximizing sequence of (\ref{eq1.8}), after scaling and translation, must
converge strongly. In Section \ref{sec4}, following Lieb we use the method of
symmetrization based on the Riesz rearrangement inequalities (\cite[section
3.7]{LiL}) and its strong form (\cite{Li1}) to show that all maximizing
functions must be radial and give another approach to the existence of
maximizing functions. In Section \ref{sec5} we use the method in \cite{Hn} to
deduce the regularity of all nonnegative critical functions. Indeed what we
will prove is a local regularity result. These results are similar in nature
to those proved in \cite{ChL,L}. In Section \ref{sec6} we use the integral
version of the method of moving planes (\cite{GNN}), which was discovered in
\cite{CLO2}, to deduce the symmetry property of the nonnegative critical
functions. Here we will need some ideas from \cite{Hn} again.

\textbf{Acknowledgment}: The research of F. Hang is supported by National
Science Foundation Grant DMS-0501050 and a Sloan Research Fellowship. The
research of X. Wang is supported by National Science Foundation Grant
DMS-0505645. The research of X. Yan is supported by National Science
Foundation Grant DMS-0401048 and an IRGP grant from Michigan State University.

\section{Basic inequalities for Poission integrals\label{sec2}}

The main aim of this section is to record some basic estimates associated with
Poission kernel and harmonic extensions which we will use freely later. For
$x_{0}\in\mathbb{R}^{n}$ and $r>0$, we write%
\[
B_{r}^{n}\left(  x_{0}\right)  =\left\{  x\in\mathbb{R}^{n}:\left\vert
x-x_{0}\right\vert <r\right\}  ,\quad B_{r}^{n}=B_{r}^{n}\left(  0\right)
,\quad B_{r}^{+}=B_{r}^{n}\cap\mathbb{R}_{+}^{n}%
\]
and $\overline{B}_{r}^{n}\left(  x_{0}\right)  $ to mean the closure of
$B_{r}^{n}\left(  x_{0}\right)  $. Assume $n\geq2$. For $t>0$, $\xi
\in\mathbb{R}^{n-1}$, we write%
\[
P_{t}\left(  \xi\right)  =\frac{2}{n\omega_{n}}\frac{t}{\left(  \left\vert
\xi\right\vert ^{2}+t^{2}\right)  ^{n/2}}.
\]
Clearly we have

\begin{itemize}
\item $P\left(  x,\xi\right)  =P_{x_{n}}\left(  x^{\prime}-\xi\right)  $ for
$x\in\mathbb{R}_{+}^{n}$, $\xi\in\mathbb{R}^{n-1}$.

\item $\left(  Pf\right)  \left(  x\right)  =\left(  P_{x_{n}}\ast f\right)
\left(  x^{\prime}\right)  $ for $x\in\mathbb{R}_{+}^{n}$.

\item $\left\vert P_{t}\right\vert _{L^{1}\left(  \mathbb{R}^{n-1}\right)
}=1$, $\left\vert P_{t}\right\vert _{L^{\infty}\left(  \mathbb{R}%
^{n-1}\right)  }=\frac{2}{n\omega_{n}}\frac{1}{t^{n-1}}$.

\item $\left\vert P_{t}\right\vert _{L^{p}\left(  \mathbb{R}^{n-1}\right)
}=c\left(  n,p\right)  t^{-\frac{\left(  n-1\right)  \left(  p-1\right)  }{p}%
}$ for $\frac{n-1}{n}<p\leq\infty$.
\end{itemize}

Recall if $X$ is a measure space, $p>0$ and $u$ is a measurable function on
$X$, then%
\[
\left\vert u\right\vert _{L_{W}^{p}\left(  X\right)  }=\sup_{t>0}t\left\vert
\left\vert u\right\vert >t\right\vert ^{1/p}.
\]
The space $L_{W}^{p}\left(  X\right)  =\left\{  u:u\text{ is measurable and
}\left\vert u\right\vert _{L_{W}^{p}\left(  X\right)  }<\infty\right\}  $.
More generally, for any $0<p<\infty$ and $0<q\leq\infty$, we have the Lorentz
norm $\left\vert \cdot\right\vert _{L^{p,q}\left(  X\right)  }$ and Lorentz
space $L^{p,q}\left(  X\right)  $ (see \cite[p188]{SW}). $L_{W}^{p}\left(
X\right)  =L^{p,\infty}\left(  X\right)  $ is a special case of such spaces.

\begin{proposition}
\label{propinterpolation}We have%
\[
\left\vert Pf\right\vert _{L_{W}^{\frac{n}{n-1}}\left(  \mathbb{R}_{+}%
^{n}\right)  }\leq c\left(  n\right)  \left\vert f\right\vert _{L^{1}\left(
\mathbb{R}^{n-1}\right)  }%
\]
and%
\begin{equation}
\left\vert Pf\right\vert _{L^{\frac{np}{n-1}}\left(  \mathbb{R}_{+}%
^{n}\right)  }\leq c\left(  n,p\right)  \left\vert f\right\vert _{L^{p}\left(
\mathbb{R}^{n-1}\right)  }\label{eq2.1}%
\end{equation}
for $1<p\leq\infty$. Moreover for $1<p<\infty$ we have%
\begin{equation}
\left\vert Pf\right\vert _{L^{\frac{np}{n-1}}\left(  \mathbb{R}_{+}%
^{n}\right)  }\leq c\left(  n,p\right)  \left\vert f\right\vert _{L^{p,\frac
{np}{n-1}}\left(  \mathbb{R}^{n-1}\right)  }\text{.}\label{eq2.2}%
\end{equation}

\begin{proof}
We only need to prove the weak type estimate. The strong estimate follows from
Marcinkiewicz interpolation theorem (see \cite[p197]{SW}) and the basic fact
$\left\vert Pf\right\vert _{L^{\infty}\left(  \mathbb{R}_{+}^{n}\right)  }%
\leq\left\vert f\right\vert _{L^{\infty}\left(  \mathbb{R}^{n-1}\right)  }$.
To prove the weak type estimate, we may assume $f\geq0$ and $\left\vert
f\right\vert _{L^{1}\left(  \mathbb{R}^{n-1}\right)  }=1$. First we observe
that $\left(  Pf\right)  \left(  x\right)  \leq\frac{c\left(  n\right)
}{x_{n}^{n-1}}$ for $x\in\mathbb{R}_{+}^{n}$ and%
\[
\int_{x\in\mathbb{R}_{+}^{n},0<x_{n}<a}\left(  Pf\right)  \left(  x\right)
dx=\int_{\mathbb{R}^{n-1}}d\xi\left(  f\left(  \xi\right)  \int_{0}^{a}%
dx_{n}\int_{\mathbb{R}^{n-1}}P\left(  x,\xi\right)  dx^{\prime}\right)  =a
\]
for $a>0$. Hence for $t>0$,%
\begin{align*}
\left\vert Pf>t\right\vert  & =\left\vert \left\{  x\in\mathbb{R}_{+}%
^{n}:0<x_{n}<c\left(  n\right)  t^{-\frac{1}{n-1}},\left(  Pf\right)  \left(
x\right)  >t\right\}  \right\vert \\
& \leq\frac{1}{t}\int_{0<x_{n}<c\left(  n\right)  t^{-\frac{1}{n-1}}%
,x^{\prime}\in\mathbb{R}^{n-1}}\left(  Pf\right)  \left(  x\right)
dx=c\left(  n\right)  t^{-\frac{n}{n-1}}.
\end{align*}
The weak type inequality follows.
\end{proof}
\end{proposition}

Later we will also need some elementary estimates for the harmonic extensions.
They are listed below.

\begin{itemize}
\item For $1\leq p\leq q\leq\infty$, we have%
\[
\left\vert P_{t}\ast f\right\vert _{L^{q}\left(  \mathbb{R}^{n-1}\right)
}\leq c\left(  n,p,q\right)  t^{-\left(  n-1\right)  \left(  \frac{1}{p}%
-\frac{1}{q}\right)  }\left\vert f\right\vert _{L^{p}\left(  \mathbb{R}%
^{n-1}\right)  }.
\]

\item Assume $f\left(  \xi\right)  =0$ for $\left\vert \xi\right\vert \geq R$,
then we have%
\[
\left\vert \left(  P_{t}\ast f\right)  \left(  \xi\right)  \right\vert
\leq\frac{c\left(  n\right)  t}{\left[  \left(  \left(  \left\vert
\xi\right\vert -R\right)  ^{+}\right)  ^{2}+t^{2}\right]  ^{n/2}}\left\vert
f\right\vert _{L^{1}\left(  \mathbb{R}^{n-1}\right)  }.
\]

\item Assume $f\left(  \xi\right)  =0$ for $\left\vert \xi\right\vert <R$,
$1\leq p\leq\infty$, then we have%
\[
\left\vert P_{t}\ast f\right\vert _{L^{\infty}\left(  B_{R/2}^{n-1}\right)
}\leq c\left(  n,p\right)  tR^{-\frac{n-1}{p}-1}\left\vert f\right\vert
_{L^{p}\left(  \mathbb{R}^{n-1}\right)  }%
\]
and%
\[
\left\vert Pf\right\vert _{L^{\infty}\left(  B_{R/2}^{+}\right)  }\leq
c\left(  n,p\right)  R^{-\frac{n-1}{p}}\left\vert f\right\vert _{L^{p}\left(
\mathbb{R}^{n-1}\right)  }.
\]

\end{itemize}

For $t>0$, $\xi\in\mathbb{R}^{n-1}$, let%
\[
Q_{t}\left(  \xi\right)  =P_{t}\left(  \xi\right)  \cdot\frac{\left\vert
\xi\right\vert }{t}=\frac{2}{n\omega_{n}}\frac{\left\vert \xi\right\vert
}{\left(  \left\vert \xi\right\vert ^{2}+t^{2}\right)  ^{n/2}},
\]
then

\begin{itemize}
\item $\left\vert Q_{t}\right\vert _{L^{p}\left(  \mathbb{R}^{n-1}\right)
}=c\left(  n,p\right)  t^{-\frac{\left(  n-1\right)  \left(  p-1\right)  }{p}%
}$ for $1<p\leq\infty$.

\item Let $\varphi\in L^{\infty}\left(  \mathbb{R}^{n-1}\right)
\cap\operatorname*{Lip}\left(  \mathbb{R}^{n-1}\right)  $, then%
\[
\left\vert P_{t}\ast\left(  \varphi f\right)  -\varphi\left(  P_{t}\ast
f\right)  \right\vert \leq\left[  \varphi\right]  _{\operatorname*{Lip}\left(
\mathbb{R}^{n-1}\right)  }tQ_{t}\ast f.
\]
In particular, it follows from Hausdorff-Young's inequality that%
\begin{align*}
& \left\vert P_{t}\ast\left(  \varphi f\right)  -\varphi\left(  P_{t}\ast
f\right)  \right\vert _{L^{q}\left(  \mathbb{R}^{n-1}\right)  }\\
& \leq c\left(  n,p,q\right)  \left[  \varphi\right]  _{\operatorname*{Lip}%
\left(  \mathbb{R}^{n-1}\right)  }t^{1-\left(  n-1\right)  \left(  \frac{1}%
{p}-\frac{1}{q}\right)  }\left\vert f\right\vert _{L^{p}\left(  \mathbb{R}%
^{n-1}\right)  }%
\end{align*}
for $1\leq p<q\leq\infty$.
\end{itemize}

As a simple application of these estimates, we derive the following
compactness result.

\begin{corollary}
For $1\leq p<\infty$, $1\leq q<\frac{np}{n-1}$, the operator%
\[
P:L^{p}\left(  \mathbb{R}^{n-1}\right)  \rightarrow L_{loc}^{q}\left(
\overline{\mathbb{R}_{+}^{n}}\right)
\]
is compact.

\begin{proof}
Assume $f_{i}\in L^{p}\left(  \mathbb{R}^{n-1}\right)  $ such that $\left\vert
f_{i}\right\vert _{L^{p}\left(  \mathbb{R}^{n-1}\right)  }\leq1$, it follows
that%
\[
\left\vert \left(  Pf_{i}\right)  \left(  x\right)  \right\vert \leq c\left(
n,p\right)  x_{n}^{-\frac{n-1}{p}}\text{ for }x\in\mathbb{R}_{+}^{n}.
\]
By the gradient estimates of harmonic functions, after passing to a
subsequence we have $Pf_{i}\rightarrow u$ in $C_{loc}^{\infty}\left(
\mathbb{R}_{+}^{n}\right)  $. For any $R>0$,%
\begin{align*}
& \left\vert Pf_{i}-Pf_{j}\right\vert _{L^{q}\left(  B_{R}^{+}\right)  }^{q}\\
& =\int_{x\in B_{R}^{+},x_{n}\geq\varepsilon}\left\vert Pf_{i}-Pf_{j}%
\right\vert ^{q}dx+\int_{x\in B_{R}^{+},x_{n}<\varepsilon}\left\vert
Pf_{i}-Pf_{j}\right\vert ^{q}dx\\
& \leq\int_{x\in B_{R}^{+},x_{n}\geq\varepsilon}\left\vert Pf_{i}%
-Pf_{j}\right\vert ^{q}dx+c\left(  n,p,q\right)  \varepsilon^{1-\left(
n-1\right)  \left(  \frac{q}{p}-1\right)  }.
\end{align*}
Hence%
\[
\lim\sup_{i,j\rightarrow\infty}\left\vert Pf_{i}-Pf_{j}\right\vert
_{L^{q}\left(  B_{R}^{+}\right)  }^{q}\leq c\left(  n,p,q\right)
\varepsilon^{1-\left(  n-1\right)  \left(  \frac{q}{p}-1\right)  }.
\]
Let $\varepsilon\rightarrow0^{+}$, we see $Pf_{i}$ is a Cauchy sequence in
$L_{loc}^{q}\left(  \overline{\mathbb{R}_{+}^{n}}\right)  $.
\end{proof}
\end{corollary}

Finally we derive a dual statement to Proposition \ref{propinterpolation}. Let
$u$ be a function on $\mathbb{R}_{+}^{n}$, we write%
\[
\left(  Tu\right)  \left(  \xi\right)  =\int_{\mathbb{R}_{+}^{n}}P\left(
x,\xi\right)  u\left(  x\right)  dx.
\]

\begin{proposition}
\label{propdualestimate}For $1\leq p<n$ we have%
\begin{equation}
\left\vert Tu\right\vert _{L^{\frac{\left(  n-1\right)  p}{n-p}}\left(
\mathbb{R}^{n-1}\right)  }\leq c\left(  n,p\right)  \left\vert u\right\vert
_{L^{p}\left(  \mathbb{R}_{+}^{n}\right)  }\label{eq2.3}%
\end{equation}
for any $u\in L^{p}\left(  \mathbb{R}_{+}^{n}\right)  $.
\end{proposition}

\begin{proof}
We may prove the inequality by a duality argument. Indeed, for any nonnegative
functions $u$ on $\mathbb{R}_{+}^{n}$ and $f$ on $\mathbb{R}^{n-1}$ we have%
\begin{align*}
0  & \leq\int_{\mathbb{R}^{n-1}}\left(  Tu\right)  \left(  \xi\right)
f\left(  \xi\right)  d\xi=\int_{\mathbb{R}^{n-1}}d\xi\int_{\mathbb{R}_{+}^{n}%
}P\left(  x,\xi\right)  u\left(  x\right)  f\left(  \xi\right)  dx\\
& =\int_{\mathbb{R}_{+}^{n}}\left(  Pf\right)  \left(  x\right)  u\left(
x\right)  dx\leq\left\vert Pf\right\vert _{L^{\frac{p}{p-1}}\left(
\mathbb{R}_{+}^{n}\right)  }\left\vert u\right\vert _{L^{p}\left(
\mathbb{R}_{+}^{n}\right)  }\\
& \leq c\left(  n,p\right)  \left\vert u\right\vert _{L^{p}\left(
\mathbb{R}_{+}^{n}\right)  }\left\vert f\right\vert _{L^{\frac{\left(
n-1\right)  p}{n\left(  p-1\right)  }}\left(  \mathbb{R}^{n-1}\right)  }.
\end{align*}
Inequality (\ref{eq2.3}) follows. We may also prove such an inequality
directly. Indeed, since%
\[
\left\vert P\left(  \cdot,\xi\right)  \right\vert _{L^{\frac{n}{n-1},\infty
}\left(  \mathbb{R}_{+}^{n}\right)  }=\left\vert P\left(  \cdot,0\right)
\right\vert _{L^{\frac{n}{n-1},\infty}\left(  \mathbb{R}_{+}^{n}\right)
}=c\left(  n\right)  <\infty,
\]
we see $T:L^{n,1}\left(  \mathbb{R}_{+}^{n}\right)  \rightarrow L^{\infty
}\left(  \mathbb{R}^{n-1}\right)  $ is a bounded linear map. On the other
hand, for $u\in L^{1}\left(  \mathbb{R}_{+}^{n}\right)  $, we have%
\[
\int_{\mathbb{R}^{n-1}}\left\vert \left(  Tu\right)  \left(  \xi\right)
\right\vert d\xi\leq\int_{\mathbb{R}^{n-1}}d\xi\int_{\mathbb{R}_{+}^{n}%
}P\left(  x,\xi\right)  \left\vert u\left(  x\right)  \right\vert
dx=\int_{\mathbb{R}_{+}^{n}}\left\vert u\left(  x\right)  \right\vert dx.
\]
Hence $T:L^{1}\left(  \mathbb{R}_{+}^{n}\right)  \rightarrow L^{1}\left(
\mathbb{R}^{n-1}\right)  $ is also bounded. The inequality (\ref{eq2.3})
follows from the Marcinkiewicz interpolation theorem.
\end{proof}

\section{The existence of maximizing functions for sharp inequalities by the
concentration compactness principle\label{sec3}}

Assume $n\geq2$ and $1<p<\infty$. Let $c_{n,p}$ be the sharp constant in
(\ref{eq2.1}), then $c_{n,p}>0$ and%
\begin{equation}
c_{n,p}^{\frac{np}{n-1}}=\sup\left\{  \int_{\mathbb{R}_{+}^{n}}\left\vert
Pf\right\vert ^{\frac{np}{n-1}}dx:f\in L^{p}\left(  \mathbb{R}^{n-1}\right)
,\left\vert f\right\vert _{L^{p}\left(  \mathbb{R}^{n-1}\right)  }=1\right\}
.\label{eq3.1}%
\end{equation}
The aim of this section is to show $c_{n,p}^{\frac{np}{n-1}}$ is attained by
some functions. Let $f$ be a function defined on $\mathbb{R}^{n-1}$. For
$\lambda>0$ and $\zeta\in\mathbb{R}^{n-1}$ we write $f^{\lambda,\zeta}\left(
\xi\right)  =\lambda^{-\frac{n-1}{p}}f\left(  \frac{\xi-\zeta}{\lambda
}\right)  $\quad for $\xi\in\mathbb{R}^{n-1}$, then%
\[
\left\vert f^{\lambda,\zeta}\right\vert _{L^{p}\left(  \mathbb{R}%
^{n-1}\right)  }=\left\vert f\right\vert _{L^{p}\left(  \mathbb{R}%
^{n-1}\right)  },\quad\left\vert Pf^{\lambda,\zeta}\right\vert _{L^{\frac
{np}{n-1}}\left(  \mathbb{R}_{+}^{n}\right)  }=\left\vert Pf\right\vert
_{L^{\frac{np}{n-1}}\left(  \mathbb{R}_{+}^{n}\right)  }.
\]
In particular the variational problem (\ref{eq3.1}) has both translation and
dilation invariance. The problem fits in the general frame of concentration
compactness principle of \cite{Lion}. We will apply this principle to prove
the following result.

\begin{theorem}
\label{thm3.1}Assume $n\geq2$ and $1<p<\infty$. Let $f_{i}$ be a maximizing
sequence of functions for (\ref{eq3.1}), then after passing to a subsequence
there exists $\lambda_{i}>0$ and $\zeta_{i}\in\mathbb{R}^{n-1} $ such that
$f_{i}^{\lambda_{i},\zeta_{i}}\rightarrow f$ in $L^{p}\left(  \mathbb{R}%
^{n-1}\right)  $. In particular, there exists at least one maximizing function
for the variational problem (\ref{eq3.1}).
\end{theorem}

A basic ingredient in the proof of Theorem \ref{thm3.1} is the following
proposition corresponding to \cite[lemma 2.1]{Lion}.

\begin{proposition}
\label{prop3.1}Assume $n\geq2$, $1<p<\infty$ and $f_{i}\in L^{p}\left(
\mathbb{R}^{n-1}\right)  $ such that $f_{i}\rightharpoonup f$ in $L^{p}\left(
\mathbb{R}^{n-1}\right)  $. After passing to a subsequence, assume%
\[
\left\vert f_{i}\right\vert ^{p}d\xi\rightharpoonup\mu\text{ in }%
\mathcal{M}\left(  \mathbb{R}^{n-1}\right)  ,\quad\left\vert Pf_{i}\right\vert
^{\frac{np}{n-1}}dx\rightharpoonup\nu\text{ in }\mathcal{M}\left(
\overline{\mathbb{R}_{+}^{n}}\right)  .
\]
Here $\mathcal{M}\left(  \mathbb{R}^{n-1}\right)  $ denotes the space of all
Radon measures on $\mathbb{R}^{n-1}$. Then we have

\begin{itemize}
\item
\[
\left.  \nu\right\vert _{\mathbb{R}_{+}^{n}}=\left\vert Pf\right\vert
^{\frac{np}{n-1}}dx.
\]
Moreover for any Borel set $E\subset\mathbb{R}^{n-1}$,%
\[
\nu\left(  E\right)  ^{\frac{n-1}{np}}\leq c_{n,p}\mu\left(  E\right)
^{\frac{1}{p}}.
\]

\item There exists a countable set of points $\zeta_{j}\in\mathbb{R}^{n-1}$
such that%
\[
\nu=\left\vert Pf\right\vert ^{\frac{np}{n-1}}dx+\sum_{j}\nu_{j}\delta
_{\zeta_{j}},\text{ }\mu\geq\left\vert f\right\vert ^{p}d\xi+\sum_{j}\mu
_{j}\delta_{\zeta_{j}},
\]
here $\mu_{j}=\mu\left(  \left\{  \zeta_{j}\right\}  \right)  $ and%
\[
\nu_{j}^{\frac{n-1}{np}}\leq c_{n,p}\mu_{j}^{\frac{1}{p}}.
\]

\item If $\nu\left(  \mathbb{R}^{n-1}\right)  ^{\frac{n-1}{np}}\geq c_{n,p}%
\mu\left(  \mathbb{R}^{n-1}\right)  ^{\frac{1}{p}}$, then $\nu$ is supported
on a single point.
\end{itemize}
\end{proposition}

\begin{proof}
Without losing of generality, we may assume $\left\vert f_{i}\right\vert
_{L^{p}\left(  \mathbb{R}^{n-1}\right)  }\leq1$. Since%
\[
\left\vert \left(  Pf_{i}\right)  \left(  x\right)  \right\vert \leq c\left(
n,p\right)  x_{n}^{-\frac{n-1}{p}}\text{\quad for }x\in\mathbb{R}_{+}^{n},
\]
it follows from the gradient estimate of harmonic functions that
$Pf_{i}\rightarrow Pf$ in $C_{loc}^{\infty}\left(  \mathbb{R}_{+}^{n}\right)
$. In particular,%
\[
\left.  \nu\right\vert _{\mathbb{R}_{+}^{n}}=\left\vert Pf\right\vert
^{\frac{np}{n-1}}dx.
\]
Let $\varphi\in C_{c}^{\infty}\left(  \mathbb{R}^{n-1}\right)  $ and $\eta\in
C_{c}^{\infty}\left(  \left[  0,\infty\right)  \right)  $ such that $0\leq
\eta\leq1$, we have%
\begin{align*}
& \left\vert \varphi\left(  x^{\prime}\right)  \eta\left(  x_{n}\right)
\left(  Pf_{i}\right)  \left(  x\right)  \right\vert _{L^{\frac{np}{n-1}%
}\left(  \mathbb{R}_{+}^{n}\right)  }\\
& \leq\left\vert \eta\left(  x_{n}\right)  P\left(  \varphi f_{i}\right)
\left(  x\right)  \right\vert _{L^{\frac{np}{n-1}}\left(  \mathbb{R}_{+}%
^{n}\right)  }+\left\vert \eta\left(  x_{n}\right)  \left(  \varphi\left(
x^{\prime}\right)  \left(  Pf_{i}\right)  \left(  x\right)  -P\left(  \varphi
f_{i}\right)  \left(  x\right)  \right)  \right\vert _{L^{\frac{np}{n-1}%
}\left(  \mathbb{R}_{+}^{n}\right)  }\\
& \leq c_{n,p}\left\vert \varphi f_{i}\right\vert _{L^{p}\left(
\mathbb{R}^{n-1}\right)  }+c\left(  n,p\right)  \left\vert \nabla
\varphi\right\vert _{L^{\infty}\left(  \mathbb{R}^{n-1}\right)  }\left(
\int_{0}^{\infty}\eta\left(  t\right)  ^{\frac{np}{n-1}}t^{\frac{np}{n-1}%
-1}dt\right)  ^{\frac{n-1}{np}}.
\end{align*}
Now fix an $\eta\in C^{\infty}\left(  \left[  0,\infty\right)  \right)  $ such
that $0\leq\eta\leq1$, $\eta\left(  0\right)  =1$ and $\eta\left(  t\right)
=0$ for $t\geq1$. For $\varepsilon>0$, denote $\eta_{\varepsilon}\left(
t\right)  =\eta\left(  t/\varepsilon\right)  $. Then%
\begin{align*}
& \left\vert \varphi\left(  x^{\prime}\right)  \eta_{\varepsilon}\left(
x_{n}\right)  \left(  Pf_{i}\right)  \left(  x\right)  \right\vert
_{L^{\frac{np}{n-1}}\left(  \mathbb{R}_{+}^{n}\right)  }\\
& \leq c_{n,p}\left\vert \varphi f_{i}\right\vert _{L^{p}\left(
\mathbb{R}^{n-1}\right)  }+c\left(  n,p\right)  \left\vert \nabla
\varphi\right\vert _{L^{\infty}\left(  \mathbb{R}^{n-1}\right)  }\varepsilon.
\end{align*}
Letting $i\rightarrow\infty$ and then $\varepsilon\rightarrow0^{+}$, we see%
\[
\left(  \int_{\mathbb{R}^{n-1}}\left\vert \varphi\right\vert ^{\frac{np}{n-1}%
}d\nu\right)  ^{\frac{n-1}{np}}\leq c_{n,p}\left(  \int_{\mathbb{R}^{n-1}%
}\left\vert \varphi\right\vert ^{p}d\mu\right)  ^{1/p}.
\]
A limit process shows for any Borel function $h$ on $\mathbb{R}^{n-1}$,%
\[
\left(  \int_{\mathbb{R}^{n-1}}\left\vert h\right\vert ^{\frac{np}{n-1}}%
d\nu\right)  ^{\frac{n-1}{np}}\leq c_{n,p}\left(  \int_{\mathbb{R}^{n-1}%
}\left\vert h\right\vert ^{p}d\mu\right)  ^{1/p}.
\]
This implies for any Borel set $E\subset\mathbb{R}^{n-1}$, $\nu\left(
E\right)  ^{\frac{n-1}{np}}\leq c_{n,p}\mu\left(  E\right)  ^{\frac{1}{p}}$.
In particular, $\nu$ is absolutely continuous with respect to $\mu$. By
Radon-Nikydom theorem (\cite[section 1.6]{EG}) we have%
\[
\nu\left(  E\right)  =\int_{E}gd\mu.
\]
Moreover for $\mu$ a.e. $\xi\in\mathbb{R}^{n-1}$%
\[
g\left(  \xi\right)  =\lim_{r\rightarrow0^{+}}\frac{\nu\left(  \overline
{B}_{r}^{n-1}\left(  \xi\right)  \right)  }{\mu\left(  \overline{B}_{r}%
^{n-1}\left(  \xi\right)  \right)  }.
\]
Let $J=\left\{  \xi\in\mathbb{R}^{n-1}:\mu\left(  \left\{  \xi\right\}
\right)  >0\right\}  $, then $J$ is countable. Moreover, for $\xi\notin J$, we
have%
\[
g\left(  \xi\right)  \leq\lim\inf_{r\rightarrow0^{+}}c_{n,p}^{\frac{np}{n-1}%
}\mu\left(  \overline{B}_{r}^{n-1}\left(  \xi\right)  \right)  ^{\frac{1}%
{n-1}}=0.
\]
Hence $\nu=\left\vert Pf\right\vert ^{\frac{np}{n-1}}dx+\sum_{j}\nu_{j}%
\delta_{\zeta_{j}}$. For the third assertion, if we know $\nu\left(
\mathbb{R}^{n-1}\right)  ^{\frac{n-1}{np}}\geq c_{n,p}\mu\left(
\mathbb{R}^{n-1}\right)  ^{\frac{1}{p}}$, then%
\[
\left(  \sum_{j}\nu_{j}\right)  ^{\frac{n-1}{np}}\geq c_{n,p}\left(  \sum
_{j}\mu_{j}\right)  ^{1/p}\geq\left(  \sum_{j}\nu_{j}^{\frac{n-1}{n}}\right)
^{1/p},
\]
hence%
\[
\left(  \sum_{j}\nu_{j}\right)  ^{\frac{n-1}{n}}\geq\sum_{j}\nu_{j}%
^{\frac{n-1}{n}}.
\]
Since $0<\frac{n-1}{n}<1$, we see at most one $\nu_{j}$ is nonzero.
\end{proof}

Now we are ready to prove Theorem \ref{thm3.1}.

\begin{proof}
[Proof of Theorem \ref{thm3.1}]For $r>0$, let%
\[
\phi_{i}\left(  r\right)  =\sup_{\zeta\in\mathbb{R}^{n-1}}\int_{B_{r}%
^{n-1}\left(  \zeta\right)  }\left\vert f_{i}\right\vert ^{p}d\xi.
\]
Then $\phi_{i}:\left(  0,\infty\right)  \rightarrow\left[  0,1\right]  $ is a
continuous nondecreasing function with%
\[
\lim_{r\rightarrow0^{+}}\phi_{i}\left(  r\right)  =0,\quad\lim_{r\rightarrow
\infty}\phi_{i}\left(  r\right)  =1.
\]
By introducing dilation factor $\lambda_{i}$ and translation by $\zeta_{i}$,
we may assume
\[
\phi_{i}\left(  1\right)  =1/2=\int_{B_{1}}\left\vert f_{i}\right\vert
^{p}d\xi.
\]
After passing to a subsequence, we may find $f\in L^{p}\left(  \mathbb{R}%
^{n-1}\right)  $ such that%
\[
f_{i}\rightharpoonup f\text{ in }L^{p}\left(  \mathbb{R}^{n-1}\right)  ,\text{
}\left\vert f_{i}\right\vert ^{p}d\xi\rightharpoonup\mu\text{ in }%
\mathcal{M}\left(  \mathbb{R}^{n-1}\right)  ,\text{ }\left\vert Pf_{i}%
\right\vert ^{\frac{np}{n-1}}dx\rightharpoonup\nu\text{ in }\mathcal{M}\left(
\overline{\mathbb{R}_{+}^{n}}\right)  .
\]
In particular, this implies $\mu\left(  \overline{B}_{1}^{n-1}\right)
\geq1/2$. We claim $\mu\left(  \mathbb{R}^{n-1}\right)  =1$. If not, then
$\mu\left(  \mathbb{R}^{n-1}\right)  =\theta\in\left(  0,1\right)  $. For
$\varepsilon>0$ small, we claim that after passing to a subsequence, we may
find $r_{0}>0$ and a sequence $r_{i}\rightarrow\infty$ such that%
\[
\theta-\varepsilon<\int_{B_{r_{0}}^{n-1}}\left\vert f_{i}\right\vert ^{p}%
d\xi\leq\int_{B_{r_{0}+2r_{i}}^{n-1}}\left\vert f_{i}\right\vert ^{p}%
d\xi<\theta+\varepsilon.
\]
Indeed, fix $r_{0}>0$ such that $\mu\left(  B_{r_{0}}^{n-1}\right)
>\theta-\varepsilon$, then for $i$ large enough, we have $\int_{B_{r_{0}%
}^{n-1}}\left\vert f_{i}\right\vert ^{p}d\xi>\theta-\varepsilon$. On the other
hand, since $\mu\left(  \overline{B}_{r_{0}+2i}^{n-1}\right)  \leq
\theta<\theta+\varepsilon$, we may inductively define $n_{i}>i$,
$n_{i+1}>n_{i}$ such that $\int_{B_{r_{0+2i}}^{n-1}}\left\vert f_{n_{i}%
}\right\vert ^{p}dx<\theta+\varepsilon$. Replacing $f_{i}$ by $f_{n_{i}}$ we
get the needed claim. Let%
\[
g_{i}=f_{i}\chi_{B_{r_{0}}^{n-1}},\quad h_{i}=f_{i}\chi_{\mathbb{R}%
^{n-1}\backslash B_{r_{0}+2r_{i}}^{n-1}}.
\]
Since%
\[
\int_{B_{r_{0}+2r_{i}}^{n-1}\backslash B_{r_{0}}^{n-1}}\left\vert
f_{i}\right\vert ^{p}d\xi\leq2\varepsilon,
\]
we see%
\[
\left\vert f_{i}-g_{i}-h_{i}\right\vert _{L^{p}\left(  \mathbb{R}%
^{n-1}\right)  }\leq c\left(  n,p\right)  \varepsilon^{1/p}.
\]
Note that%
\begin{align*}
& \left\vert \left\vert Pg_{i}+Ph_{i}\right\vert ^{\frac{np}{n-1}}-\left\vert
Pg_{i}\right\vert ^{\frac{np}{n-1}}-\left\vert Ph_{i}\right\vert ^{\frac
{np}{n-1}}\right\vert \\
& \leq c\left(  n,p\right)  \left(  \left\vert Pg_{i}\right\vert ^{\frac
{np}{n-1}-1}\left\vert Ph_{i}\right\vert +\left\vert Pg_{i}\right\vert
\left\vert Ph_{i}\right\vert ^{\frac{np}{n-1}-1}\right)  .
\end{align*}
On the other hand,%
\begin{align*}
& \int_{\mathbb{R}_{+}^{n}}\left\vert Pg_{i}\right\vert ^{\frac{np}{n-1}%
-1}\left\vert Ph_{i}\right\vert dx\\
& \leq\left\vert Pg_{i}\right\vert _{L^{\frac{np}{n-1}}\left(  B_{R}%
^{+}\right)  }^{\frac{np}{n-1}-1}\left\vert Ph_{i}\right\vert _{L^{\frac
{np}{n-1}}\left(  B_{R}^{+}\right)  }+\left\vert Pg_{i}\right\vert
_{L^{\frac{np}{n-1}}\left(  \mathbb{R}_{+}^{n}\backslash B_{R}^{+}\right)
}^{\frac{np}{n-1}-1}\left\vert Ph_{i}\right\vert _{L^{\frac{np}{n-1}}\left(
\mathbb{R}_{+}^{n}\backslash B_{R}^{+}\right)  }\\
& \leq c\left(  n,p\right)  R^{\frac{n-1}{p}}r_{i}^{-\frac{n-1}{p}}+c\left(
n,p\right)  r_{0}^{\frac{\left(  p-1\right)  \left(  np-n+1\right)  }{p}%
}\left\vert \frac{x_{n}}{\left[  \left(  \left(  \left\vert x^{\prime
}\right\vert -r_{0}\right)  ^{+}\right)  ^{2}+x_{n}^{2}\right]  ^{n/2}%
}\right\vert _{L^{\frac{np}{n-1}}\left(  \mathbb{R}_{+}^{n}\backslash
B_{R}^{+}\right)  }^{\frac{np}{n-1}-1},
\end{align*}
this implies%
\begin{align*}
& \lim\sup_{i\rightarrow\infty}\int_{\mathbb{R}_{+}^{n}}\left\vert
Pg_{i}\right\vert ^{\frac{np}{n-1}-1}\left\vert Ph_{i}\right\vert dx\\
& \leq c\left(  n,p\right)  r_{0}^{\frac{\left(  p-1\right)  \left(
np-n+1\right)  }{p}}\left\vert \frac{x_{n}}{\left[  \left(  \left(  \left\vert
x^{\prime}\right\vert -r_{0}\right)  ^{+}\right)  ^{2}+x_{n}^{2}\right]
^{n/2}}\right\vert _{L^{\frac{np}{n-1}}\left(  \mathbb{R}_{+}^{n}\backslash
B_{R}^{+}\right)  }^{\frac{np}{n-1}-1}.
\end{align*}
Let $R\rightarrow\infty$, we see%
\[
\lim_{i\rightarrow\infty}\int_{\mathbb{R}_{+}^{n}}\left\vert Pg_{i}\right\vert
^{\frac{np}{n-1}-1}\left\vert Ph_{i}\right\vert dx=0.
\]
Similarly,%
\[
\lim_{i\rightarrow\infty}\int_{\mathbb{R}_{+}^{n}}\left\vert Pg_{i}\right\vert
\left\vert Ph_{i}\right\vert ^{\frac{np}{n-1}-1}dx=0.
\]
Hence%
\[
\lim_{i\rightarrow\infty}\int_{\mathbb{R}_{+}^{n}}\left\vert \left\vert
Pg_{i}+Ph_{i}\right\vert ^{\frac{np}{n-1}}-\left\vert Pg_{i}\right\vert
^{\frac{np}{n-1}}-\left\vert Ph_{i}\right\vert ^{\frac{np}{n-1}}\right\vert
dx=0.
\]
Since%
\begin{align*}
\int_{\mathbb{R}_{+}^{n}}\left\vert Pg_{i}\right\vert ^{\frac{np}{n-1}}dx  &
\leq c_{n,p}^{\frac{np}{n-1}}\left\vert g_{i}\right\vert _{L^{p}\left(
\mathbb{R}^{n-1}\right)  }^{\frac{np}{n-1}}\leq c_{n,p}^{\frac{np}{n-1}%
}\left(  \theta+\varepsilon\right)  ^{\frac{n}{n-1}},\\
\int_{\mathbb{R}_{+}^{n}}\left\vert Ph_{i}\right\vert ^{\frac{np}{n-1}}dx  &
\leq c_{n,p}^{\frac{np}{n-1}}\left\vert h_{i}\right\vert _{L^{p}\left(
\mathbb{R}^{n-1}\right)  }^{\frac{np}{n-1}}\leq c_{n,p}^{\frac{np}{n-1}%
}\left(  1-\theta+\varepsilon\right)  ^{\frac{n}{n-1}},
\end{align*}
we see%
\begin{align*}
& c_{n,p}^{\frac{np}{n-1}}+o\left(  1\right) \\
& =\int_{\mathbb{R}_{+}^{n}}\left\vert Pf_{i}\right\vert ^{\frac{np}{n-1}}dx\\
& \leq\left(  \left\vert Pg_{i}+Ph_{i}\right\vert _{L^{\frac{np}{n-1}}\left(
\mathbb{R}_{+}^{n}\right)  }+c\left(  n,p\right)  \varepsilon^{1/p}\right)
^{\frac{np}{n-1}}\\
& \leq\int_{\mathbb{R}_{+}^{n}}\left\vert Pg_{i}+Ph_{i}\right\vert ^{\frac
{np}{n-1}}dx+c\left(  n,p\right)  \varepsilon^{1/p}\\
& \leq\int_{\mathbb{R}_{+}^{n}}\left(  \left\vert Pg_{i}\right\vert
^{\frac{np}{n-1}}+\left\vert Ph_{i}\right\vert ^{\frac{np}{n-1}}\right)
dx+c\left(  n,p\right)  \varepsilon^{1/p}+o\left(  1\right) \\
& \leq c_{n,p}^{\frac{np}{n-1}}\left(  \theta+\varepsilon\right)  ^{\frac
{n}{n-1}}+c_{n,p}^{\frac{np}{n-1}}\left(  1-\theta+\varepsilon\right)
^{\frac{n}{n-1}}+c\left(  n,p\right)  \varepsilon^{1/p}+o\left(  1\right)  .
\end{align*}
Letting $i\rightarrow\infty$ and then $\varepsilon\rightarrow0^{+}$, we see%
\[
1\leq\theta^{\frac{n}{n-1}}+\left(  1-\theta\right)  ^{\frac{n}{n-1}}.
\]
This gives us a contradiction since $\frac{n}{n-1}>1$. Hence $\mu\left(
\mathbb{R}^{n-1}\right)  =1$. Next we claim $\nu\left(  \overline
{\mathbb{R}_{+}^{n}}\right)  =c_{n,p}^{\frac{np}{n-1}}$. Indeed, for any
$\varepsilon>0$ small, we may find $r>0$ such that $\mu\left(  B_{r}%
^{n-1}\right)  >1-\varepsilon$, this implies $\int_{B_{r}^{n-1}}\left\vert
f_{i}\right\vert ^{p}d\xi>1-\varepsilon$ when $i$ is large enough. Hence
$\int_{\mathbb{R}^{n-1}\backslash B_{r}^{n-1}}\left\vert f_{i}\right\vert
^{p}d\xi\leq\varepsilon$. Let $g_{i}=f_{i}\chi_{B_{r}^{n-1}}$ and $h_{i}%
=f_{i}\chi_{\mathbb{R}^{n-1}\backslash B_{r}^{n-1}}$, then%
\begin{align*}
\left\vert Ph_{i}\right\vert _{L^{\frac{np}{n-1}}\left(  \mathbb{R}_{+}%
^{n}\right)  }  & \leq c\left(  n,p\right)  \varepsilon^{1/p},\\
\left\vert \left(  Pg_{i}\right)  \left(  x\right)  \right\vert  & \leq
c\left(  n,p\right)  r^{\frac{\left(  n-1\right)  \left(  p-1\right)  }{p}%
}\frac{x_{n}}{\left[  \left(  \left(  \left\vert x^{\prime}\right\vert
-r\right)  ^{+}\right)  ^{2}+x_{n}^{2}\right]  ^{n/2}}.
\end{align*}
This implies%
\begin{align*}
& \int_{\mathbb{R}^{n}\backslash B_{R}^{+}}\left\vert Pf_{i}\right\vert
^{\frac{np}{n-1}}dx\\
& \leq c\left(  n,p\right)  \varepsilon^{\frac{n}{n-1}}+c\left(  n,p\right)
r^{n\left(  p-1\right)  }\left\vert \frac{x_{n}}{\left[  \left(  \left(
\left\vert x^{\prime}\right\vert -r\right)  ^{+}\right)  ^{2}+x_{n}%
^{2}\right]  ^{n/2}}\right\vert _{L^{\frac{np}{n-1}}\left(  \mathbb{R}%
^{n}\backslash B_{R}^{+}\right)  }^{\frac{np}{n-1}}.
\end{align*}
Taking a limit for $i\rightarrow\infty$, we see%
\begin{align*}
\nu\left(  \overline{\mathbb{R}_{+}^{n}}\right)   & \geq\nu\left(
\overline{B_{R}^{+}}\right) \\
& \geq c_{n,p}^{\frac{np}{n-1}}-c\left(  n,p\right)  \varepsilon^{\frac
{n}{n-1}}-c\left(  n,p\right)  r^{n\left(  p-1\right)  }\left\vert \frac
{x_{n}}{\left[  \left(  \left(  \left\vert x^{\prime}\right\vert -r\right)
^{+}\right)  ^{2}+x_{n}^{2}\right]  ^{n/2}}\right\vert _{L^{\frac{np}{n-1}%
}\left(  \mathbb{R}^{n}\backslash B_{R}^{+}\right)  }^{\frac{np}{n-1}}.
\end{align*}
let $R\rightarrow\infty$ then $\varepsilon\rightarrow0^{+}$, we see
$\nu\left(  \overline{\mathbb{R}_{+}^{n}}\right)  =c_{n,p}^{\frac{np}{n-1}}$.

By Proposition \ref{prop3.1} we know there exists a countable set of points
$\zeta_{j}\in\mathbb{R}^{n-1}$ such that%
\[
\nu=\left\vert Pf\right\vert ^{\frac{np}{n-1}}dx+\sum_{j}\nu_{j}\delta
_{\zeta_{j}},\text{\quad}\mu\geq\left\vert f\right\vert ^{p}d\xi+\sum_{j}%
\mu_{j}\delta_{\zeta_{j}},
\]
here $\mu_{j}=\mu\left(  \left\{  \zeta_{j}\right\}  \right)  $ and%
\[
\nu_{j}^{\frac{n-1}{np}}\leq c_{n,p}\mu_{j}^{\frac{1}{p}}.
\]
If $f=0$, then $\nu\left(  \mathbb{R}^{n-1}\right)  =c_{n,p}^{\frac{np}{n-1}}$
and hence $\nu\left(  \mathbb{R}^{n-1}\right)  ^{\frac{n-1}{np}}=c_{n,p}%
\mu\left(  \mathbb{R}^{n-1}\right)  ^{\frac{1}{p}}$. This implies for some
$\zeta_{1}\in\mathbb{R}^{n-1}$, $\nu=c_{n,p}^{\frac{np}{n-1}}\delta_{\zeta
_{1}} $. In particular, $\mu\left(  \left\{  \zeta_{1}\right\}  \right)  \geq1
$ and this implies $\mu=\delta_{\zeta_{1}}$. But%
\[
\int_{B_{1}^{n-1}\left(  \zeta_{1}\right)  }\left\vert f_{i}\right\vert
^{p}d\xi\leq1/2
\]
implies $\mu\left(  B_{1}^{n-1}\left(  \zeta_{1}\right)  \right)  \leq1/2$.
This gives us a contradiction. Hence $f\neq0$. Now%
\[
c_{n,p}^{\frac{np}{n-1}}=\left\vert Pf\right\vert _{L^{\frac{np}{n-1}}\left(
\mathbb{R}_{+}^{n}\right)  }^{\frac{np}{n-1}}+\sum_{j}\nu_{j}\leq
c_{n,p}^{\frac{np}{n-1}}\left\vert f\right\vert _{L^{p}\left(  \mathbb{R}%
^{n-1}\right)  }^{\frac{np}{n-1}}+c_{n,p}^{\frac{np}{n-1}}\sum_{j}\mu
_{j}^{\frac{n}{n-1}},
\]
hence%
\[
1\leq\left\vert f\right\vert _{L^{p}\left(  \mathbb{R}^{n-1}\right)  }%
^{\frac{np}{n-1}}+\sum_{j}\mu_{j}^{\frac{n}{n-1}}.
\]
But since%
\[
1\geq\left\vert f\right\vert _{L^{p}\left(  \mathbb{R}^{n-1}\right)  }%
^{p}+\sum_{j}\mu_{j}%
\]
and $\frac{n}{n-1}>1$, we see $\mu_{j}=0$ and $\left\vert f\right\vert
_{L^{p}\left(  \mathbb{R}^{n-1}\right)  }=1$. This implies $f_{i}\rightarrow
f$ in $L^{p}\left(  \mathbb{R}^{n-1}\right)  $.
\end{proof}

\section{The existence of maximizing functions for sharp inequalities by
symmetrization\label{sec4}}

Following Lieb (\cite{Li2}), using the method of symmetrization we will show
all the maximizers of variational problem (\ref{eq3.1}) are radial symmetric
with respect to some point and we will give another approach to the existence
of maximizing functions.

Let $u$ be a measurable function on $\mathbb{R}^{n}$, the symmetric
rearrangement of $u$ is the nonnegative lower semi-continuous radial
decreasing function $u^{\ast}$ which has the same distribution as $u$. It
satisfies the following important Riesz rearrangement inequality
(\cite[p87]{LiL}): for any nonnegative measurable functions $u,v,w$ on
$\mathbb{R}^{n}$, we have%
\[
\int_{\mathbb{R}^{n}}dx\int_{\mathbb{R}^{n}}u\left(  x\right)  v\left(
y-x\right)  w\left(  y\right)  dy\leq\int_{\mathbb{R}^{n}}dx\int
_{\mathbb{R}^{n}}u^{\ast}\left(  x\right)  v^{\ast}\left(  y-x\right)
w^{\ast}\left(  y\right)  dy.
\]
Using the fact $\left\vert w\right\vert _{L^{p}\left(  \mathbb{R}^{n}\right)
}=\left\vert w^{\ast}\right\vert _{L^{p}\left(  \mathbb{R}^{n}\right)  }$ for
$p>0$, we see for $1\leq p\leq\infty$,%
\[
\left\vert u\ast v\right\vert _{L^{p}\left(  \mathbb{R}^{n}\right)  }%
\leq\left\vert u^{\ast}\ast v^{\ast}\right\vert _{L^{p}\left(  \mathbb{R}%
^{n}\right)  }.
\]
Moreover if $u$ is nonnegative radial symmetric and strictly decreasing in the
radial direction, $v$ is nonnegative, $1<p<\infty$ and
\[
\left\vert u\ast v\right\vert _{L^{p}\left(  \mathbb{R}^{n}\right)
}=\left\vert u\ast v^{\ast}\right\vert _{L^{p}\left(  \mathbb{R}^{n}\right)
}<\infty,
\]
then for some $x_{0}\in\mathbb{R}^{n}$, we have $v\left(  x\right)  =v^{\ast
}\left(  x-x_{0}\right)  $.

Indeed, we may assume $v$ is not identically zero. Choose a nonnegative $w\in
L^{p^{\prime}}\left(  \mathbb{R}^{n}\right)  $ with $\left\vert w\right\vert
_{L^{p^{\prime}}\left(  \mathbb{R}^{n}\right)  }=1$ such that%
\[
\left\vert u\ast v\right\vert _{L^{p}\left(  \mathbb{R}^{n}\right)  }%
=\int_{\mathbb{R}^{n}}\left(  u\ast v\right)  \left(  y\right)  w\left(
y\right)  dy.
\]
Then we have%
\begin{align*}
\left\vert u\ast v\right\vert _{L^{p}\left(  \mathbb{R}^{n}\right)  }  &
=\int_{\mathbb{R}^{n}}dx\int_{\mathbb{R}^{n}}u\left(  x\right)  v\left(
y-x\right)  w\left(  y\right)  dy\\
& \leq\int_{\mathbb{R}^{n}}dx\int_{\mathbb{R}^{n}}u\left(  x\right)  v^{\ast
}\left(  y-x\right)  w^{\ast}\left(  y\right)  dy\\
& =\int_{\mathbb{R}^{n}}\left(  u\ast v^{\ast}\right)  \left(  y\right)
w^{\ast}\left(  y\right)  dy\\
& \leq\left\vert u\ast v^{\ast}\right\vert _{L^{p}\left(  \mathbb{R}%
^{n}\right)  }=\left\vert u\ast v\right\vert _{L^{p}\left(  \mathbb{R}%
^{n}\right)  },
\end{align*}
hence all the inequalities become equalities. It follows from the Lieb's
strong version of Riesz rearrangement inequality (\cite{Li1}) that for some
$x_{0}\in\mathbb{R}^{n}$, $v\left(  x\right)  =v^{\ast}\left(  x-x_{0}\right)
$.

\begin{theorem}
\label{thm4.1}Assume $n\geq2$ and $1<p<\infty$, then the value%
\[
c_{n,p}^{\frac{np}{n-1}}=\sup\left\{  \int_{\mathbb{R}_{+}^{n}}\left\vert
Pf\right\vert ^{\frac{np}{n-1}}dx:f\in L^{p}\left(  \mathbb{R}^{n-1}\right)
,\left\vert f\right\vert _{L^{p}\left(  \mathbb{R}^{n-1}\right)  }=1\right\}
,
\]
is attained by some functions. After multiplying by a nonzero constant, every
maximizer $f$ is nonnegative, radial symmetric with respect to some point,
strictly decreasing in the radial direction and it satisfies%
\begin{equation}
f\left(  \xi\right)  ^{p-1}=\int_{\mathbb{R}_{+}^{n}}P\left(  x,\xi\right)
\left(  Pf\right)  \left(  x\right)  ^{\frac{np}{n-1}-1}dx.\label{eq4.1}%
\end{equation}

If $n\geq3$ and $p=\frac{2\left(  n-1\right)  }{n-2}$, then any maximizer must
be of the form%
\[
f\left(  \xi\right)  =\pm c\left(  n\right)  \left(  \frac{\lambda}%
{\lambda^{2}+\left\vert \xi-\xi_{0}\right\vert ^{2}}\right)  ^{\frac{n-2}{2}}%
\]
for some $\lambda>0$, $\xi_{0}\in\mathbb{R}^{n-1}$. In particular,
$c_{n,\frac{2\left(  n-1\right)  }{n-2}}=n^{-\frac{n-2}{2\left(  n-1\right)
}}\omega_{n}^{-\frac{n-2}{2n\left(  n-1\right)  }}$.

If $n\geq3$ and $p=\frac{2\left(  n-1\right)  }{n}$, then any maximizer must
be of the form%
\[
f\left(  \xi\right)  =\pm c\left(  n\right)  \left(  \frac{\lambda}%
{\lambda^{2}+\left\vert \xi-\xi_{0}\right\vert ^{2}}\right)  ^{\frac{n}{2}}%
\]
for some $\lambda>0$, $\xi_{0}\in\mathbb{R}^{n-1}$. In particular,
$c_{n,\frac{2\left(  n-1\right)  }{n}}=\frac{1}{\sqrt{2\left(  n-2\right)
}\sqrt[4]{\pi}}\left(  \frac{\left(  n-2\right)  !}{\Gamma\left(  \frac
{n-1}{2}\right)  }\right)  ^{\frac{1}{2\left(  n-1\right)  }}$.
\end{theorem}

\begin{proof}
Assume $f_{i}$ is a maximizing sequence. Since $\left\vert f_{i}^{\ast
}\right\vert _{L^{p}\left(  \mathbb{R}^{n-1}\right)  }=\left\vert
f_{i}\right\vert _{L^{p}\left(  \mathbb{R}^{n-1}\right)  }=1$ and
\begin{align*}
\left\vert Pf_{i}\right\vert _{L^{\frac{np}{n-1}}\left(  \mathbb{R}_{+}%
^{n}\right)  }^{\frac{np}{n-1}}  & =\int_{0}^{\infty}\left\vert P_{x_{n}}\ast
f_{i}\right\vert _{L^{\frac{np}{n-1}}\left(  \mathbb{R}^{n-1}\right)  }%
^{\frac{np}{n-1}}dx_{n}\\
& \leq\int_{0}^{\infty}\left\vert P_{x_{n}}\ast f_{i}^{\ast}\right\vert
_{L^{\frac{np}{n-1}}\left(  \mathbb{R}^{n-1}\right)  }^{\frac{np}{n-1}}%
dx_{n}=\left\vert Pf_{i}^{\ast}\right\vert _{L^{\frac{np}{n-1}}\left(
\mathbb{R}_{+}^{n}\right)  }^{\frac{np}{n-1}},
\end{align*}
we see $f_{i}^{\ast}$ is again a maximizing sequence. Hence we may assume
$f_{i}$ is a nonnegative radial decreasing function.

For any $f\in L^{p}\left(  \mathbb{R}^{n-1}\right)  $ and any $\lambda>0$, we
let $f^{\lambda}\left(  \xi\right)  =\lambda^{-\frac{n-1}{p}}f\left(
\frac{\xi}{\lambda}\right)  $, then it is clear that $\left\vert f^{\lambda
}\right\vert _{L^{p}\left(  \mathbb{R}^{n-1}\right)  }=\left\vert f\right\vert
_{L^{p}\left(  \mathbb{R}^{n-1}\right)  }$ and $\left\vert Pf^{\lambda
}\right\vert _{L^{\frac{np}{n-1}}\left(  \mathbb{R}_{+}^{n}\right)
}=\left\vert Pf\right\vert _{L^{\frac{np}{n-1}}\left(  \mathbb{R}_{+}%
^{n}\right)  }$. For convenience, denote $e_{1}=\left(  1,0,\cdots,0\right)
\in\mathbb{R}^{n}$ and%
\[
a_{i}=\sup_{\lambda>0}f_{i}^{\lambda}\left(  e_{1}^{\prime}\right)
=\sup_{\lambda>0}\lambda^{-\frac{n-1}{p}}f_{i}\left(  \frac{e_{1}^{\prime}%
}{\lambda}\right)  .
\]
It follows that%
\[
0\leq f_{i}\left(  \xi\right)  \leq a_{i}\left\vert \xi\right\vert
^{-\frac{n-1}{p}}%
\]
and hence%
\[
\left\vert f_{i}\right\vert _{L^{p,\infty}\left(  \mathbb{R}^{n-1}\right)
}\leq\omega_{n-1}^{1/p}a_{i}.
\]
Now%
\begin{align*}
\left\vert Pf_{i}\right\vert _{L^{\frac{np}{n-1}}\left(  \mathbb{R}_{+}%
^{n}\right)  }  & \leq c\left(  n,p\right)  \left\vert f_{i}\right\vert
_{L^{p,\frac{np}{n-1}}\left(  \mathbb{R}_{+}^{n}\right)  }\\
& \leq c\left(  n,p\right)  \left\vert f_{i}\right\vert _{L^{p}\left(
\mathbb{R}_{+}^{n}\right)  }^{\frac{n-1}{n}}\left\vert f_{i}\right\vert
_{L^{p,\infty}\left(  \mathbb{R}_{+}^{n}\right)  }^{\frac{1}{n}}\\
& \leq c\left(  n,p\right)  a_{i}^{1/n},
\end{align*}
this implies $a_{i}\geq c\left(  n,p\right)  >0$. We may choose $\lambda_{i}>0
$ such that $f_{i}^{\lambda_{i}}\left(  e_{1}^{\prime}\right)  \geq c\left(
n,p\right)  >0$. Replacing $f_{i}$ by $f_{i}^{\lambda_{i}}$ we may assume
$f\left(  e_{1}^{\prime}\right)  \geq c\left(  n,p\right)  >0$. On the other
hand, since $f_{i}$ is nonnegative radial decreasing and $\left\vert
f_{i}\right\vert _{L^{p}\left(  \mathbb{R}^{n-1}\right)  }=1$, we see%
\[
\left\vert f_{i}\left(  \xi\right)  \right\vert \leq\omega_{n-1}%
^{-1/p}\left\vert \xi\right\vert ^{-\left(  n-1\right)  /p}.
\]
Hence after passing to a subsequence, we may find a nonnegative radial
decreasing function $f$ such that $f_{i}\rightarrow f$ a.e.. It follows that
$f\left(  \xi\right)  \geq c\left(  n,p\right)  >0$ for $\left\vert
\xi\right\vert \leq1$, $f_{i}\rightharpoonup f$ in $L^{p}\left(
\mathbb{R}^{n-1}\right)  $ and $\left\vert f\right\vert _{L^{p}\left(
\mathbb{R}^{n-1}\right)  }\leq1$. Since%
\[
\int_{\mathbb{R}^{n-1}}\left\vert \left\vert f_{i}\right\vert ^{p}-\left\vert
f\right\vert ^{p}-\left\vert f_{i}-f\right\vert ^{p}\right\vert d\xi
\rightarrow0,
\]
we see%
\begin{align*}
\left\vert f_{i}-f\right\vert _{L^{p}\left(  \mathbb{R}^{n-1}\right)  }^{p}  &
=\left\vert f_{i}\right\vert _{L^{p}\left(  \mathbb{R}^{n-1}\right)  }%
^{p}-\left\vert f\right\vert _{L^{p}\left(  \mathbb{R}^{n-1}\right)  }%
^{p}+o\left(  1\right) \\
& =1-\left\vert f\right\vert _{L^{p}\left(  \mathbb{R}^{n-1}\right)  }%
^{p}+o\left(  1\right)  .
\end{align*}
On the other hand, since $\left(  Pf_{i}\right)  \left(  x\right)
\rightarrow\left(  Pf\right)  \left(  x\right)  $ for $x\in\mathbb{R}_{+}^{n}$
and $\left\vert Pf_{i}\right\vert _{L^{\frac{np}{n-1}}\left(  \mathbb{R}%
_{+}^{n}\right)  }\leq c_{n,p}$, we see%
\begin{align*}
\left\vert Pf_{i}\right\vert _{L^{\frac{np}{n-1}}\left(  \mathbb{R}_{+}%
^{n}\right)  }^{\frac{np}{n-1}}  & =\left\vert Pf\right\vert _{L^{\frac
{np}{n-1}}\left(  \mathbb{R}_{+}^{n}\right)  }^{\frac{np}{n-1}}+\left\vert
Pf_{i}-Pf\right\vert _{L^{\frac{np}{n-1}}\left(  \mathbb{R}_{+}^{n}\right)
}^{\frac{np}{n-1}}+o\left(  1\right) \\
& \leq c_{n,p}^{\frac{np}{n-1}}\left\vert f\right\vert _{L^{p}\left(
\mathbb{R}^{n-1}\right)  }^{\frac{np}{n-1}}+c_{n,p}^{\frac{np}{n-1}}\left\vert
f_{i}-f\right\vert _{L^{p}\left(  \mathbb{R}^{n-1}\right)  }^{\frac{np}{n-1}%
}+o\left(  1\right)  .
\end{align*}
Hence%
\[
1\leq\left\vert f\right\vert _{L^{p}\left(  \mathbb{R}^{n-1}\right)  }%
^{\frac{np}{n-1}}+\left\vert f_{i}-f\right\vert _{L^{p}\left(  \mathbb{R}%
^{n-1}\right)  }^{\frac{np}{n-1}}+o\left(  1\right)  .
\]
Let $i\rightarrow\infty$, we see
\[
1\leq\left\vert f\right\vert _{L^{p}\left(  \mathbb{R}^{n-1}\right)  }%
^{\frac{np}{n-1}}+\left(  1-\left\vert f\right\vert _{L^{p}\left(
\mathbb{R}^{n-1}\right)  }^{p}\right)  ^{\frac{n}{n-1}}.
\]
Since $\frac{n}{n-1}>1$ and $f\neq0$, we see $\left\vert f\right\vert
_{L^{p}\left(  \mathbb{R}^{n-1}\right)  }=1$. Hence $f_{i}\rightarrow f$ in
$L^{p}\left(  \mathbb{R}^{n-1}\right)  $ and $f$ is a maximizer. This implies
the existence of an extremal function.

Assume $f\in L^{p}\left(  \mathbb{R}^{n-1}\right)  $ is a maximizer, then so
is $\left\vert f\right\vert $. Hence $\left\vert Pf\right\vert _{L^{\frac
{np}{n-1}}\left(  \mathbb{R}_{+}^{n}\right)  }=\left\vert P\left\vert
f\right\vert \right\vert _{L^{\frac{np}{n-1}}\left(  \mathbb{R}_{+}%
^{n}\right)  }$. On the other hand, since $\left\vert \left(  Pf\right)
\left(  x\right)  \right\vert \leq P\left(  \left\vert f\right\vert \right)
\left(  x\right)  $ for $x\in\mathbb{R}_{+}^{n}$, we see $\left\vert
Pf\right\vert =P\left(  \left\vert f\right\vert \right)  $ and this implies
either $f\geq0$ or $f\leq0$. Assume $f\geq0$, then the Euler-Lagrange equation
is given by%
\[
\int_{\mathbb{R}_{+}^{n}}P\left(  x,\xi\right)  \left(  Pf\right)  \left(
x\right)  ^{\frac{np}{n-1}-1}dx=c\cdot f\left(  \xi\right)  ^{p-1}.
\]
Here $c$ is a constant. Using the fact $\left\vert f\right\vert _{L^{p}\left(
\mathbf{R}^{n-1}\right)  }=1$, we see%
\[
c=\left\vert Pf\right\vert _{L^{\frac{np}{n-1}}\left(  \mathbb{R}_{+}%
^{n}\right)  }^{\frac{np}{n-1}}=c_{n,p}^{\frac{np}{n-1}}.
\]
After scaling by a positive constant we get%
\[
\int_{\mathbb{R}_{+}^{n}}P\left(  x,\xi\right)  \left(  Pf\right)  \left(
x\right)  ^{\frac{np}{n-1}-1}dx=f\left(  \xi\right)  ^{p-1}.
\]
On the other hand, we know for $x_{n}>0$, $\left\vert P_{x_{n}}\ast
f\right\vert _{L^{\frac{np}{n-1}}\left(  \mathbb{R}^{n-1}\right)  }=\left\vert
P_{x_{n}}\ast f^{\ast}\right\vert _{L^{\frac{np}{n-1}}\left(  \mathbb{R}%
^{n-1}\right)  }$, this implies $f\left(  \xi\right)  =f^{\ast}\left(  \xi
-\xi_{0}\right)  $ for some $\xi_{0}$. It follows from the Euler-Lagrange
equation that $f$ must be strictly decreasing along the radial direction.

For the case when $p=\frac{2\left(  n-1\right)  }{n-2}$, we first observe that
if $f\in L^{\frac{2\left(  n-1\right)  }{n-2}}\left(  \mathbb{R}^{n-1}\right)
$, let $u=Pf$, $\widetilde{f}\left(  \xi\right)  =\frac{1}{\left\vert
\xi\right\vert ^{n-2}}f\left(  \frac{\xi}{\left\vert \xi\right\vert ^{2}%
}\right)  $ and $\widetilde{u}\left(  x\right)  =\frac{1}{\left\vert
x\right\vert ^{n-2}}u\left(  \frac{x}{\left\vert x\right\vert ^{2}}\right)  $,
then we have $\widetilde{u}=P\widetilde{f}$, $\left\vert \widetilde
{f}\right\vert _{L^{\frac{2\left(  n-1\right)  }{n-2}}\left(  \mathbb{R}%
^{n-1}\right)  }=\left\vert f\right\vert _{L^{\frac{2\left(  n-1\right)
}{n-2}}\left(  \mathbb{R}^{n-1}\right)  }$ and $\left\vert \widetilde
{u}\right\vert _{L^{\frac{2n}{n-2}}\left(  \mathbb{R}_{+}^{n}\right)
}=\left\vert u\right\vert _{L^{\frac{2n}{n-2}}\left(  \mathbb{R}_{+}%
^{n}\right)  }$. This is the conformal invariance property for the particular
power. As a consequence, if $f$ is a maximizer which is nonnegative and
radial, then $\frac{1}{\left\vert \xi\right\vert ^{n-2}}f\left(  \frac{\xi
}{\left\vert \xi\right\vert ^{2}}-e_{1}^{\prime}\right)  $ is a maximizer too.
In particular, $\frac{1}{\left\vert \xi\right\vert ^{n-2}}f\left(  \frac{\xi
}{\left\vert \xi\right\vert ^{2}}-e_{1}^{\prime}\right)  $ is radial with
respect to some point. To find such $f$, we prove the following facts.
\end{proof}

\begin{proposition}
\label{prop4.1}Let $n\geq2$, $u$ be a function on $\mathbb{R}^{n}$ which is
radial with respect to the origin, $0<u\left(  x\right)  <\infty$ for $x\neq0
$, $e_{1}=\left(  1,0,\cdots,0\right)  $, $\alpha\in\mathbb{R}$, $\alpha\neq
0$. If $v\left(  x\right)  =\left\vert x\right\vert ^{\alpha}u\left(  \frac
{x}{\left\vert x\right\vert ^{2}}-e_{1}\right)  $ is radial with respect to
some point, then either $u\left(  x\right)  =\left(  c_{1}\left\vert
x\right\vert ^{2}+c_{2}\right)  ^{\alpha/2}$ for some $c_{1}\geq0$, $c_{2}>0 $
or%
\[
u\left(  x\right)  =\left\{
\begin{tabular}
[c]{l}%
$c_{1}\left\vert x\right\vert ^{\alpha},$ if $x\neq0,$\\
$c_{2},$ if $x=0,$%
\end{tabular}
\right.
\]
for some $c_{1}>0$ and $c_{2}$, an arbitrary number.

\begin{proof}
First we observe that $\left\vert \frac{x}{\left\vert x\right\vert ^{2}}%
-e_{1}\right\vert =1$ if and only if $x_{1}=\frac{1}{2}$. For $r>0$, $r\neq1
$, we have $\left\vert \frac{x}{\left\vert x\right\vert ^{2}}-e_{1}\right\vert
=r$ if and only if $x\in\partial B_{\frac{r}{\left\vert r^{2}-1\right\vert }%
}\left(  \frac{e_{1}}{1-r^{2}}\right)  $. By scaling, we may assume $u\left(
e_{1}\right)  =1$. Then $v\left(  \frac{1}{2},x^{\prime\prime}\right)
=\left(  \frac{1}{4}+\left\vert x^{\prime\prime}\right\vert ^{2}\right)
^{\alpha/2}$. Assume $v$ is symmetric with respect to $z=\left(
z_{1},z^{\prime\prime}\right)  $. Then $v\left(  \frac{1}{2},\cdot\right)  $
is symmetric with respect to $z^{\prime\prime}$, hence $z^{\prime\prime}=0$.
Denote $z=ae_{1}$, we claim $0\leq a\leq1$. If this is not the case, then we
may find a $r>0$, $r\neq1$ such that $a=\frac{1}{1-r^{2}}$. Now on $\partial
B_{\frac{r}{\left\vert r^{2}-1\right\vert }}\left(  \frac{e_{1}}{1-r^{2}%
}\right)  $, $v\left(  x\right)  =\left\vert x\right\vert ^{\alpha}u\left(
re_{1}\right)  $ and it is not a constant function, contradiction. For
$x=\left(  \frac{1}{2},x^{\prime\prime}\right)  $, we have%
\[
v\left(  x\right)  =\left(  \left\vert x-ae_{1}\right\vert ^{2}+a-a^{2}%
\right)  ^{\alpha/2}.
\]
Hence%
\[
v\left(  x\right)  =\left(  \left\vert x-ae_{1}\right\vert ^{2}+a-a^{2}%
\right)  ^{\alpha/2}=\left(  \left\vert x\right\vert ^{2}-2ax_{1}+a\right)
^{\alpha/2}%
\]
for $\left\vert x-ae_{1}\right\vert \geq\left\vert \frac{1}{2}-a\right\vert $.
When $a=\frac{1}{2}$, we see $v\left(  x\right)  =\left(  \left\vert
x\right\vert ^{2}-2ax_{1}+a\right)  ^{\alpha/2}$ for all $x$. This implies
$u\left(  x\right)  =\left(  \frac{1}{2}\left\vert x\right\vert ^{2}+\frac
{1}{2}\right)  ^{\alpha/2}$. Hence we assume $a\neq\frac{1}{2}$ from now on.
Without losing of generality, we assume $0\leq a<\frac{1}{2}$. We claim that
\begin{equation}
v\left(  x\right)  =\left(  \left\vert x\right\vert ^{2}-2ax_{1}+a\right)
^{\alpha/2}\label{eq4.2}%
\end{equation}
for all $x\neq0$. To see this, we first make the following observation. Assume
for some given $r>0$, $r\neq1$ and for some $y\in\partial B_{\frac
{r}{\left\vert r^{2}-1\right\vert }}\left(  \frac{e_{1}}{1-r^{2}}\right)  $,
(\ref{eq4.2}) is true for $y$, then it is true for all $x\in\partial
B_{\frac{r}{\left\vert r^{2}-1\right\vert }}\left(  \frac{e_{1}}{1-r^{2}%
}\right)  $. Indeed, for $x$ on such a sphere, we have
\[
\frac{1-2x_{1}}{\left\vert x\right\vert ^{2}}=r^{2}-1.
\]
Hence
\begin{align*}
v\left(  x\right)   & =\left\vert x\right\vert ^{\alpha}u\left(
re_{1}\right)  =\left\vert x\right\vert ^{\alpha}\left\vert y\right\vert
^{-\alpha}\left(  \left\vert y\right\vert ^{2}-2ay_{1}+a\right)  ^{\alpha/2}\\
& =\left\vert x\right\vert ^{\alpha}\left(  1+\frac{a\left(  1-2x_{1}\right)
}{\left\vert x\right\vert ^{2}}\right)  ^{\alpha/2}=\left(  \left\vert
x\right\vert ^{2}-2ax_{1}+a\right)  ^{\alpha/2}.
\end{align*}
Note that we know (\ref{eq4.2}) is true for $x=te_{1}$ with $t\in\left(
-\infty,-\frac{1}{2}\right]  $. By the above observation we know it is also
true for $te_{1}$ with $t\in\left[  \frac{1}{4},\frac{1}{2}\right]  $. This
implies it is true for $te_{1}$ with $t\in\left(  -\infty,-\frac{1}{4}\right]
$. Go back we see it is true for $te_{1}$ with $t\in\left[  \frac{1}{6}%
,\frac{1}{2}\right]  $. Keep this procedure going, we see (\ref{eq4.2}) is
true for all $te_{1}$ with $t\neq0$. Hence it is true for all $x\neq0$. This
implies $u\left(  x\right)  =\left(  a\left\vert x\right\vert ^{2}+1-a\right)
^{\alpha/2}$.
\end{proof}
\end{proposition}

\begin{remark}
\label{rmk4.1}The case when $\alpha=0$ is a little bit different. However one
has: Let $n\geq2$, $u$ be a function on $\mathbb{R}^{n}$ which is radial with
respect to the origin, $e_{1}=\left(  1,0,\cdots,0\right)  $. If $v\left(
x\right)  =u\left(  \frac{x}{\left\vert x\right\vert ^{2}}-e_{1}\right)  $ is
radial with respect to some point, then either%
\[
u\left(  x\right)  =\left\{
\begin{tabular}
[c]{l}%
$c_{1},$ if $x=0,$\\
$c_{2},$ if $x\neq0,$%
\end{tabular}
\right.
\]
or there exists $r>0$, $r\neq1$ such that%
\[
u\left(  x\right)  =\left\{
\begin{array}
[c]{l}%
c_{1},\text{ if }\left\vert x\right\vert <r,\\
c_{2},\text{ if }\left\vert x\right\vert =r,\\
c_{3},\text{ if }\left\vert x\right\vert >r.
\end{array}
\right.
\]
Here $c_{i}$'s are arbitrary constants.
\end{remark}

\begin{proof}
[Proof of Theorem \ref{thm4.1} continued]Since $\left\vert f\right\vert
_{L^{\frac{2\left(  n-1\right)  }{n-2}}\left(  \mathbb{R}^{n-1}\right)  }=1$
and it is strictly decreasing along the radial direction, we see $0<f\left(
\xi\right)  <\infty$ for $\xi\neq0$. Note that since $f$ satisfies the
Euler-Lagrange equation, it is defined everywhere instead of almost
everywhere. It follows from Proposition \ref{prop4.1} that $f\left(
\xi\right)  =\left(  c_{1}\left\vert \xi\right\vert ^{2}+c_{2}\right)
^{-\frac{n-2}{2}}$ for some $c_{1},c_{2}>0$ (note that $f$ can not be a
constant function and the scalar multiple of $\left\vert \xi\right\vert
^{2-n}$ is ruled out by the integrability). A simple change of variable shows%
\[
1=\int_{\mathbb{R}^{n-1}}f\left(  \xi\right)  ^{\frac{2\left(  n-1\right)
}{n-2}}d\xi=\left(  c_{1}c_{2}\right)  ^{-\frac{n-1}{2}}\int_{\mathbb{R}%
^{n-1}}\left(  1+\left\vert \xi\right\vert ^{2}\right)  ^{-\left(  n-1\right)
}d\xi.
\]
Hence $c_{1}c_{2}=c\left(  n\right)  $. It follows that for some $\lambda>0$,
$f\left(  \xi\right)  =c\left(  n\right)  \left(  \frac{\lambda}{\lambda
^{2}+\left\vert \xi\right\vert ^{2}}\right)  ^{\frac{n-2}{2}}$. Let
$e_{n}=\left(  0,\cdots,0,1\right)  $. Since $u\left(  x\right)  =\left\vert
x+e_{n}\right\vert ^{2-n}$ is a bounded harmonic function on $\overline
{\mathbb{R}_{+}^{n}}$ and $u\left(  \xi,0\right)  =\left(  1+\left\vert
\xi\right\vert ^{2}\right)  ^{-\frac{n-2}{2}}$, we see%
\[
P\left(  \left(  1+\left\vert \xi\right\vert ^{2}\right)  ^{-\frac{n-2}{2}%
}\right)  \left(  x\right)  =\left\vert x+e_{n}\right\vert ^{2-n}.
\]
By the dilation invariance%
\[
c_{n,\frac{2\left(  n-1\right)  }{n-2}}=\frac{\left\vert \left\vert
x+e_{n}\right\vert ^{2-n}\right\vert _{L^{\frac{2n}{n-2}}\left(
\mathbb{R}_{+}^{n}\right)  }}{\left\vert \left(  1+\left\vert \xi\right\vert
^{2}\right)  ^{-\frac{n-2}{2}}\right\vert _{L^{\frac{2\left(  n-1\right)
}{n-2}}\left(  \mathbb{R}^{n-1}\right)  }}=n^{-\frac{n-2}{2\left(  n-1\right)
}}\omega_{n}^{-\frac{n-2}{2n\left(  n-1\right)  }}.
\]

For the case when $p=\frac{2\left(  n-1\right)  }{n}$, we know any maximizer
after multiplying by a constant will be nonnegative and satisfy%
\[
f\left(  \xi\right)  ^{\frac{n-2}{n}}=\int_{\mathbb{R}_{+}^{n}}P\left(
x,\xi\right)  \left(  Pf\right)  \left(  x\right)  dx=c\left(  n\right)
\int_{\mathbb{R}^{n-1}}\frac{f\left(  \zeta\right)  }{\left\vert \xi
-\zeta\right\vert ^{n-2}}d\zeta.
\]
Let $g\left(  \xi\right)  =f\left(  \xi\right)  ^{\frac{n-2}{n}}$, then $g\in
L^{\frac{2\left(  n-1\right)  }{n-2}}\left(  \mathbb{R}^{n-1}\right)  $ and%
\[
g\left(  \xi\right)  =c\left(  n\right)  \int_{\mathbb{R}^{n-1}}\frac{g\left(
\zeta\right)  ^{\frac{\left(  n-1\right)  +1}{\left(  n-1\right)  -1}}%
}{\left\vert \xi-\zeta\right\vert ^{\left(  n-1\right)  -1}}d\zeta.
\]
It follows from \cite[theorem 1]{CLO2} or \cite{L} that for some $\lambda>0$
and $\xi_{0}\in\mathbb{R}^{n-1}$, we have%
\[
g\left(  \xi\right)  =c\left(  n\right)  \left(  \frac{\lambda}{\lambda
^{2}+\left\vert \xi-\xi_{0}\right\vert ^{2}}\right)  ^{\frac{n-2}{2}}.
\]
Hence%
\[
f\left(  \xi\right)  =c\left(  n\right)  \left(  \frac{\lambda}{\lambda
^{2}+\left\vert \xi-\xi_{0}\right\vert ^{2}}\right)  ^{\frac{n}{2}}.
\]
Since $u\left(  x\right)  =\frac{x_{n}+1}{\left\vert x+e_{n}\right\vert ^{n}}$
is a bounded harmonic function on $\overline{\mathbb{R}_{+}^{n}}$ and
$u\left(  \xi,0\right)  =\left(  1+\left\vert \xi\right\vert ^{2}\right)
^{-\frac{n}{2}}$, we see%
\[
P\left(  \left(  1+\left\vert \xi\right\vert ^{2}\right)  ^{-\frac{n}{2}%
}\right)  \left(  x\right)  =\frac{x_{n}+1}{\left\vert x+e_{n}\right\vert
^{n}}.
\]
By the dilation invariance%
\[
c_{n,\frac{2\left(  n-1\right)  }{n}}=\frac{\left\vert \frac{x_{n}%
+1}{\left\vert x+e_{n}\right\vert ^{n}}\right\vert _{L^{2}\left(
\mathbb{R}_{+}^{n}\right)  }}{\left\vert \left(  1+\left\vert \xi\right\vert
^{2}\right)  ^{-\frac{n}{2}}\right\vert _{L^{\frac{2\left(  n-1\right)  }{n}%
}\left(  \mathbb{R}^{n-1}\right)  }}=\frac{1}{\sqrt{2\left(  n-2\right)
}\sqrt[4]{\pi}}\left(  \frac{\left(  n-2\right)  !}{\Gamma\left(  \frac
{n-1}{2}\right)  }\right)  ^{\frac{1}{2\left(  n-1\right)  }}.
\]

\end{proof}

As a final note, we point out the similar statement to Proposition
\ref{prop4.1} in dimension one.

\begin{proposition}
\label{prop4.2}Assume $u\in C^{3}\left(  \mathbb{R}\right)  $, $u>0$,
$\alpha\in\mathbb{R}$ such that for any $y\in\mathbb{R}$, $\left\vert
x\right\vert ^{\alpha}u\left(  \frac{1}{x}+y\right)  $ is symmetric with
respect to some point, then for some $a\geq0$, $b>0$ and $x_{0}\in\mathbb{R}$,
we have $u\left(  x\right)  =\left[  a\left(  x-x_{0}\right)  ^{2}+b\right]
^{\alpha/2}$.

\begin{proof}
Assume $\left\vert x\right\vert ^{\alpha}u\left(  \frac{1}{x}+y\right)  $ is
symmetric with respect to $z=z\left(  y\right)  $, then%
\[
\left\vert x\right\vert ^{\alpha}u\left(  \frac{1}{x}+y\right)  =\left\vert
2z-x\right\vert ^{\alpha}u\left(  \frac{1}{2z-x}+y\right)  .
\]
Replace $x$ by $x^{-1}$, we see%
\[
u\left(  x+y\right)  =\left\vert 1-2zx\right\vert ^{\alpha}u\left(  y-\frac
{x}{1-2zx}\right)  .
\]
Calculation shows%
\begin{align*}
& \left\vert 1-2zx\right\vert ^{\alpha}u\left(  y-\frac{x}{1-2zx}\right) \\
& =u\left(  y\right)  -\left(  u^{\prime}\left(  y\right)  +2\alpha zu\left(
y\right)  \right)  x+\left(  \frac{u^{\prime\prime}\left(  y\right)  }%
{2}+2\left(  \alpha-1\right)  zu^{\prime}\left(  y\right)  +2\alpha\left(
\alpha-1\right)  z^{2}u\left(  y\right)  \right)  x^{2}\\
& -\left[  \frac{u^{\prime\prime\prime}\left(  y\right)  }{6}+\left(
\alpha-2\right)  zu^{\prime\prime}\left(  y\right)  +2\left(  \alpha-1\right)
\left(  \alpha-2\right)  z^{2}u^{\prime}\left(  y\right)  +\frac{4}{3}%
\alpha\left(  \alpha-1\right)  \left(  \alpha-2\right)  z^{3}u\left(
y\right)  \right]  x^{3}\\
& +o\left(  x^{3}\right)  .
\end{align*}
Comparing the Taylor expansion coefficients, we see%
\[
u^{\prime}\left(  y\right)  =-\alpha zu\left(  y\right)
\]
and%
\[
\frac{u^{\prime\prime\prime}\left(  y\right)  }{3}+\left(  \alpha-2\right)
zu^{\prime\prime}\left(  y\right)  +2\left(  \alpha-1\right)  \left(
\alpha-2\right)  z^{2}u^{\prime}\left(  y\right)  +\frac{4}{3}\alpha\left(
\alpha-1\right)  \left(  \alpha-2\right)  z^{3}u\left(  y\right)  =0
\]
If $\alpha=0$, then we see $u^{\prime}=0$ and hence $u$ must be a constant
function and we are done. Assume $\alpha\neq0$, then%
\[
z=-\frac{u^{\prime}\left(  y\right)  }{\alpha u\left(  y\right)  }.
\]
Plug this in the second equation, we get%
\[
u^{2}u^{\prime\prime\prime}+3\left(  \frac{2}{\alpha}-1\right)  u{}u^{\prime
}u^{\prime\prime}+\left(  \frac{2}{\alpha}-1\right)  \left(  \frac{2}{\alpha
}-2\right)  u^{\prime3}=0.
\]
Hence%
\[
\left(  u^{2/\alpha}\right)  ^{\prime\prime\prime}=\frac{2}{\alpha}u^{\frac
{2}{\alpha}-3}\left[  u^{2}u^{\prime\prime\prime}+3\left(  \frac{2}{\alpha
}-1\right)  u{}u^{\prime}u^{\prime\prime}+\left(  \frac{2}{\alpha}-1\right)
\left(  \frac{2}{\alpha}-2\right)  u^{\prime3}\right]  =0.
\]
The proposition follows.
\end{proof}
\end{proposition}

\section{Regularity of nonnegative critical functions\label{sec5}}

In this section we will study the regularity issue related to the
Euler-Lagrange equation (\ref{eq1.9}). Let $f$ be a nonnegative function
satisfying (\ref{eq1.9}), define $u=Pf$, then the single equation becomes an
integral system%
\begin{align*}
u\left(  x\right)   & =\int_{\mathbb{R}^{n-1}}P\left(  x,\xi\right)  f\left(
\xi\right)  d\xi,\\
f\left(  \xi\right)  ^{p-1}  & =\int_{\mathbb{R}_{+}^{n}}P\left(
x,\xi\right)  u\left(  x\right)  ^{\frac{np}{n-1}-1}dx.
\end{align*}
This system is very similar to the one appeared in the study of the sharp
Hardy-Littlewood-Sobolev inequality (\cite[part (ii) of theorem 2.3]{Li2}). In
\cite{ChL,L} the regularity problem for some cases of that system was resolved
by a linear approach. In \cite{Hn}, a nonlinear approach was introduced to
resolve the regularity issue for all the cases. We will apply the nonlinear
approach to handle (\ref{eq1.9}).

\begin{theorem}
\label{thm5.1}Assume $n\geq2$, $1<p<\infty$, $f\in L_{loc}^{p}\left(
\mathbb{R}^{n-1}\right)  $ is nonnegative, not identically zero and it
satisfies%
\[
f\left(  \xi\right)  ^{p-1}=\int_{\mathbb{R}_{+}^{n}}P\left(  x,\xi\right)
\left(  Pf\right)  \left(  x\right)  ^{\frac{np}{n-1}-1}dx,
\]
then $f\in C^{\infty}\left(  \mathbb{R}^{n-1}\right)  $. If we know $f\in
L^{p}\left(  \mathbb{R}^{n-1}\right)  $, then $f\left(  \xi\right)
\rightarrow0$ as $\left\vert \xi\right\vert \rightarrow\infty$.
\end{theorem}

We note that the condition $f\in L_{loc}^{p}\left(  \mathbb{R}^{n-1}\right)  $
can not be dropped, since the above equation has $c\left(  n,p\right)
\left\vert \xi\right\vert ^{-\frac{n-1}{p}}$ as a singular solution. To prove
this theorem, we first derive some local regularity results for some integral
inequalities. According to the range of $p$, we need two local results stated
in Proposition \ref{prop5.1} and Proposition \ref{prop5.2} below. The two
propositions are of the same nature as \cite[proposition 2.1]{Hn} and
\cite[theorem 1.3]{L}.

\begin{proposition}
\label{prop5.1}Given $n\geq2$, $1<a,b\leq\infty$, $1\leq r<\infty$, $\frac
{n}{n-1}<p<q<\infty$ such that%
\[
\frac{1}{n}<\frac{r}{q}+\frac{1}{a}<\frac{r}{p}+\frac{1}{a}\leq1
\]
and%
\[
\frac{n}{ra}+\frac{n-1}{b}=\frac{1}{r}.
\]
Denote $B_{R}=B_{R}^{n-1}$ and $B_{R}^{+}=B_{R}^{n}\cap\mathbb{R}_{+}^{n}$.
Assume $u,v\in L^{p}\left(  B_{R}^{+}\right)  $, $U\in L^{a}\left(  B_{R}%
^{+}\right)  $, $F\in L^{b}\left(  B_{R}\right)  $ are all nonnegative
functions with $\left.  v\right\vert _{B_{R/2}^{+}}\in L^{q}\left(
B_{R/2}^{+}\right)  $,%
\[
\left\vert U\right\vert _{L^{a}\left(  B_{R}^{+}\right)  }^{1/r}\left\vert
F\right\vert _{L^{b}\left(  B_{R}\right)  }\leq\varepsilon\left(
n,p,q,r,a,b\right)  \text{ small}%
\]
and%
\[
u\left(  x\right)  \leq\int_{B_{R}}P\left(  x,\xi\right)  F\left(  \xi\right)
\left[  \int_{B_{R}^{+}}P\left(  y,\xi\right)  U\left(  y\right)  u\left(
y\right)  ^{r}dy\right]  ^{1/r}d\xi+v\left(  x\right)
\]
for $x\in B_{R}^{+}$, then we have $\left.  u\right\vert _{B_{R/4}^{+}}\in
L^{q}\left(  B_{R/4}^{+}\right)  $ and%
\[
\left\vert u\right\vert _{L^{q}\left(  B_{R/4}^{+}\right)  }\leq c\left(
n,p,q,r,a,b\right)  \left(  R^{\frac{n}{q}-\frac{n}{p}}\left\vert u\right\vert
_{L^{p}\left(  B_{R}^{+}\right)  }+\left\vert v\right\vert _{L^{q}\left(
B_{R/2}^{+}\right)  }\right)  .
\]

\end{proposition}

\begin{proof}
By scaling we may assume $R=1$. First assume we have $u,v\in L^{q}\left(
B_{1}^{+}\right)  $. Denote%
\[
f\left(  \xi\right)  =\int_{B_{1}^{+}}P\left(  x,\xi\right)  U\left(
x\right)  u\left(  x\right)  ^{r}dx\text{ for }\xi\in B_{1}.
\]
Let $p_{1}$ and $q_{1}$ be the numbers defined by%
\[
\frac{n-1}{p_{1}}=\frac{nr}{p}+\frac{n}{a}-1,\quad\frac{n-1}{q_{1}}=\frac
{nr}{q}+\frac{n}{a}-1,
\]
then it follows from Proposition \ref{propdualestimate} that%
\begin{align*}
\left\vert f\right\vert _{L^{p_{1}}\left(  B_{1}\right)  }  & \leq c\left(
n,p,r,a\right)  \left\vert U\right\vert _{L^{a}\left(  B_{1}^{+}\right)
}\left\vert u\right\vert _{L^{p}\left(  B_{1}^{+}\right)  }^{r},\\
\left\vert f\right\vert _{L^{q_{1}}\left(  B_{1}\right)  }  & \leq c\left(
n,q,r,a\right)  \left\vert U\right\vert _{L^{a}\left(  B_{1}^{+}\right)
}\left\vert u\right\vert _{L^{q}\left(  B_{1}^{+}\right)  }^{r}.
\end{align*}
Given $0<s<t\leq1/2$. For $x\in B_{s}^{+}$, we have%
\begin{align*}
& u\left(  x\right) \\
& \leq\int_{B_{\frac{s+t}{2}}}P\left(  x,\xi\right)  F\left(  \xi\right)
f\left(  \xi\right)  ^{1/r}d\xi+\int_{B_{1}\backslash B_{\frac{s+t}{2}}%
}P\left(  x,\xi\right)  F\left(  \xi\right)  f\left(  \xi\right)  ^{1/r}%
d\xi+v\left(  x\right) \\
& \leq\int_{B_{\frac{s+t}{2}}}P\left(  x,\xi\right)  F\left(  \xi\right)
f\left(  \xi\right)  ^{1/r}d\xi+\frac{c\left(  n\right)  }{\left(  t-s\right)
^{n-1}}\int_{B_{1}\backslash B_{\frac{s+t}{2}}}F\left(  \xi\right)  f\left(
\xi\right)  ^{1/r}d\xi+v\left(  x\right) \\
& \leq\int_{B_{\frac{s+t}{2}}}P\left(  x,\xi\right)  F\left(  \xi\right)
f\left(  \xi\right)  ^{1/r}d\xi+\frac{c\left(  n,p\right)  }{\left(
t-s\right)  ^{n-1}}\left\vert F\right\vert _{L^{b}\left(  B_{1}\right)
}\left\vert f\right\vert _{L^{p_{1}}\left(  B_{1}\right)  }^{1/r}+v\left(
x\right) \\
& \leq\int_{B_{\frac{s+t}{2}}}P\left(  x,\xi\right)  F\left(  \xi\right)
f\left(  \xi\right)  ^{1/r}d\xi+\frac{c\left(  n,p,r,a\right)  }{\left(
t-s\right)  ^{n-1}}\left\vert u\right\vert _{L^{p}\left(  B_{1}^{+}\right)
}+v\left(  x\right)  .
\end{align*}
Hence%
\[
\left\vert u\right\vert _{L^{q}\left(  B_{s}^{+}\right)  }\leq c\left(
n,q\right)  \left\vert F\right\vert _{L^{b}\left(  B_{1}\right)  }\left\vert
f\right\vert _{L^{q_{1}}\left(  B_{\frac{s+t}{2}}\right)  }^{1/r}%
+\frac{c\left(  n,p,q,r,a\right)  }{\left(  t-s\right)  ^{n-1}}\left\vert
u\right\vert _{L^{p}\left(  B_{1}^{+}\right)  }+\left\vert v\right\vert
_{L^{q}\left(  B_{1/2}^{+}\right)  }.
\]
On the other hand, for $\xi\in B_{\frac{s+t}{2}}$, we have%
\begin{align*}
f\left(  \xi\right)   & =\int_{B_{t}^{+}}P\left(  x,\xi\right)  U\left(
x\right)  u\left(  x\right)  ^{r}dx+\int_{B_{1}^{+}\backslash B_{t}^{+}%
}P\left(  x,\xi\right)  U\left(  x\right)  u\left(  x\right)  ^{r}dx\\
& \leq\int_{B_{t}^{+}}P\left(  x,\xi\right)  U\left(  x\right)  u\left(
x\right)  ^{r}dx+\frac{c\left(  n\right)  }{\left(  t-s\right)  ^{n-1}}%
\int_{B_{1}^{+}\backslash B_{t}^{+}}U\left(  x\right)  u\left(  x\right)
^{r}dx\\
& \leq\int_{B_{t}^{+}}P\left(  x,\xi\right)  U\left(  x\right)  u\left(
x\right)  ^{r}dx+\frac{c\left(  n,p,r,a\right)  }{\left(  t-s\right)  ^{n-1}%
}\left\vert U\right\vert _{L^{a}\left(  B_{1}^{+}\right)  }\left\vert
u\right\vert _{L^{p}\left(  B_{1}^{+}\right)  }^{r}.
\end{align*}
It follows from Proposition \ref{propdualestimate} that%
\[
\left\vert f\right\vert _{L^{q_{1}}\left(  B_{\frac{s+t}{2}}\right)  }\leq
c\left(  n,q,r,a\right)  \left\vert U\right\vert _{L^{a}\left(  B_{1}%
^{+}\right)  }\left\vert u\right\vert _{L^{q}\left(  B_{t}^{+}\right)  }%
^{r}+\frac{c\left(  n,p,q,r,a\right)  }{\left(  t-s\right)  ^{n-1}}\left\vert
U\right\vert _{L^{a}\left(  B_{1}^{+}\right)  }\left\vert u\right\vert
_{L^{p}\left(  B_{1}^{+}\right)  }^{r}.
\]
Combining the two inequalities together, we see%
\[
\left\vert u\right\vert _{L^{q}\left(  B_{s}^{+}\right)  }\leq\frac{1}%
{2}\left\vert u\right\vert _{L^{q}\left(  B_{t}^{+}\right)  }+\frac{c\left(
n,p,q,r,a\right)  }{\left(  t-s\right)  ^{n-1}}\left\vert u\right\vert
_{L^{p}\left(  B_{1}^{+}\right)  }+\left\vert v\right\vert _{L^{q}\left(
B_{1/2}^{+}\right)  }%
\]
if $\varepsilon$ is small enough. It follows from the usual iteration
procedure (\cite[lemma 4.3 on p.75]{HL}) that%
\[
\left\vert u\right\vert _{L^{q}\left(  B_{1/4}^{+}\right)  }\leq c\left(
n,p,q,r,a\right)  \left(  \left\vert u\right\vert _{L^{p}\left(  B_{1}%
^{+}\right)  }+\left\vert v\right\vert _{L^{q}\left(  B_{1/2}^{+}\right)
}\right)  .
\]
To prove the full proposition, we note that there exists a function $0\leq
\eta\left(  x\right)  \leq1$ such that%
\[
u\left(  x\right)  =\eta\left(  x\right)  \int_{B_{1}}P\left(  x,\xi\right)
F\left(  \xi\right)  \left[  \int_{B_{1}^{+}}P\left(  y,\xi\right)  U\left(
y\right)  u\left(  y\right)  ^{r}dy\right]  ^{1/r}d\xi+\eta\left(  x\right)
v\left(  x\right)  .
\]
We may define a map $T$ by%
\[
T\left(  \varphi\right)  \left(  x\right)  =\eta\left(  x\right)  \int_{B_{1}%
}P\left(  x,\xi\right)  F\left(  \xi\right)  \left[  \int_{B_{1}^{+}}P\left(
y,\xi\right)  U\left(  y\right)  \left\vert \varphi\left(  y\right)
\right\vert ^{r}dy\right]  ^{1/r}d\xi.
\]
Note that we have%
\begin{align*}
\left\vert T\left(  \varphi\right)  \right\vert _{L^{p}\left(  B_{1}%
^{+}\right)  }  & \leq c\left(  n,p,r,a,b\right)  \left\vert U\right\vert
_{L^{a}\left(  B_{1}^{+}\right)  }^{1/r}\left\vert F\right\vert _{L^{b}\left(
B_{1}\right)  }\left\vert \varphi\right\vert _{L^{p}\left(  B_{1}\right)
}\leq\frac{1}{2}\left\vert \varphi\right\vert _{L^{p}\left(  B_{1}\right)
},\\
\left\vert T\left(  \varphi\right)  \right\vert _{L^{q}\left(  B_{1}%
^{+}\right)  }  & \leq c\left(  n,q,r,a,b\right)  \left\vert U\right\vert
_{L^{a}\left(  B_{1}^{+}\right)  }^{1/r}\left\vert F\right\vert _{L^{b}\left(
B_{1}\right)  }\left\vert \varphi\right\vert _{L^{q}\left(  B_{1}\right)
}\leq\frac{1}{2}\left\vert \varphi\right\vert _{L^{q}\left(  B_{1}\right)  }%
\end{align*}
if $\varepsilon$ is small enough. Moreover, for $\varphi,\psi\in L^{p}\left(
B_{1}^{+}\right)  $, it follows from Minkowski inequality that%
\[
\left\vert T\left(  \varphi\right)  \left(  x\right)  -T\left(  \psi\right)
\left(  x\right)  \right\vert \leq T\left(  \left\vert \varphi-\psi\right\vert
\right)  \left(  x\right)  \text{ for }x\in B_{1}^{+},
\]
hence%
\[
\left\vert T\left(  \varphi\right)  -T\left(  \psi\right)  \right\vert
_{L^{p}\left(  B_{1}^{+}\right)  }\leq\left\vert T\left(  \left\vert
\varphi-\psi\right\vert \right)  \right\vert _{L^{p}\left(  B_{1}^{+}\right)
}\leq\frac{1}{2}\left\vert \varphi-\psi\right\vert _{L^{p}\left(  B_{1}%
^{+}\right)  }.
\]
Similarly we have for any $\varphi,\psi\in L^{q}\left(  B_{1}^{+}\right)  $,%
\[
\left\vert T\left(  \varphi\right)  -T\left(  \psi\right)  \right\vert
_{L^{q}\left(  B_{1}^{+}\right)  }\leq\frac{1}{2}\left\vert \varphi
-\psi\right\vert _{L^{q}\left(  B_{1}^{+}\right)  }.
\]
For $k\in\mathbb{N}$, let $v_{k}\left(  x\right)  =\min\left\{  v\left(
x\right)  ,k\right\}  $, then it follows from contraction mapping theorem that
we may find a unique $u_{k}\in L^{q}\left(  B_{1}^{+}\right)  $ such that%
\begin{align*}
& u_{k}\left(  x\right) \\
& =T\left(  u_{k}\right)  \left(  x\right)  +\eta\left(  x\right)
v_{k}\left(  x\right) \\
& =\eta\left(  x\right)  \int_{B_{1}}P\left(  x,\xi\right)  F\left(
\xi\right)  \left[  \int_{B_{1}^{+}}P\left(  y,\xi\right)  U\left(  y\right)
\left\vert u_{k}\left(  y\right)  \right\vert ^{r}dy\right]  ^{1/r}d\xi
+\eta\left(  x\right)  v_{k}\left(  x\right)  .
\end{align*}
Applying the apriori estimate to $u_{k}$, we see%
\[
\left\vert u_{k}\right\vert _{L^{q}\left(  B_{1/4}^{+}\right)  }\leq c\left(
n,p,q,r,a\right)  \left(  \left\vert u_{k}\right\vert _{L^{p}\left(  B_{1}%
^{+}\right)  }+\left\vert v\right\vert _{L^{q}\left(  B_{1/2}^{+}\right)
}\right)  .
\]
Observe that%
\[
u\left(  x\right)  =T\left(  u\right)  \left(  x\right)  +\eta\left(
x\right)  v\left(  x\right)  ,
\]
we see%
\begin{align*}
\left\vert u_{k}-u\right\vert _{L^{p}\left(  B_{1}^{+}\right)  }  &
\leq\left\vert T\left(  u_{k}\right)  -T\left(  u\right)  \right\vert
_{L^{p}\left(  B_{1}^{+}\right)  }+\left\vert v_{k}-v\right\vert
_{L^{p}\left(  B_{1}^{+}\right)  }\\
& \leq\frac{1}{2}\left\vert u_{k}-u\right\vert _{L^{p}\left(  B_{1}%
^{+}\right)  }+\left\vert v_{k}-v\right\vert _{L^{p}\left(  B_{1}^{+}\right)
}.
\end{align*}
Hence $\left\vert u_{k}-u\right\vert _{L^{p}\left(  B_{1}^{+}\right)  }%
\leq2\left\vert v_{k}-v\right\vert _{L^{p}\left(  B_{1}^{+}\right)
}\rightarrow0$ as $k\rightarrow\infty$. Take a limit process in the apriori
estimate, we get the proposition.
\end{proof}

The other local regularity result is

\begin{proposition}
\label{prop5.2}Given $n\geq2$, $1<a,b\leq\infty$, $1\leq r<\infty$,
$1<p<q<\infty$ such that%
\[
0<\frac{r}{p}+\frac{1}{a}<1
\]
and%
\[
\frac{n-1}{ra}+\frac{n}{b}=1.
\]
Denote $B_{R}=B_{R}^{n-1}$ and $B_{R}^{+}=B_{R}^{n}\cap\mathbb{R}_{+}^{n}$.
Assume $f,g\in L^{p}\left(  B_{R}\right)  $, $F\in L^{a}\left(  B_{R}\right)
$, $U\in L^{b}\left(  B_{R}^{+}\right)  $ are all nonnegative functions with
$\left.  g\right\vert _{B_{R/2}}\in L^{q}\left(  B_{R/2}\right)  $,%
\[
\left\vert F\right\vert _{L^{a}\left(  B_{R}\right)  }^{1/r}\left\vert
U\right\vert _{L^{b}\left(  B_{R}^{+}\right)  }\leq\varepsilon\left(
n,p,q,r,a,b\right)  \text{ small}%
\]
and%
\[
f\left(  \xi\right)  \leq\int_{B_{R}^{+}}P\left(  x,\xi\right)  U\left(
x\right)  \left[  \int_{B_{R}}P\left(  x,\zeta\right)  F\left(  \zeta\right)
f\left(  \zeta\right)  ^{r}d\zeta\right]  ^{1/r}dx+g\left(  \xi\right)
\]
for $\xi\in B_{R}$, then we have $\left.  f\right\vert _{B_{R/4}}\in
L^{q}\left(  B_{R/4}\right)  $ and%
\[
\left\vert f\right\vert _{L^{q}\left(  B_{R/4}\right)  }\leq c\left(
n,p,q,r,a,b\right)  \left(  R^{\frac{n-1}{q}-\frac{n-1}{p}}\left\vert
f\right\vert _{L^{p}\left(  B_{R}\right)  }+\left\vert g\right\vert
_{L^{q}\left(  B_{R/2}\right)  }\right)  .
\]

\end{proposition}

\begin{proof}
By scaling, we may assume $R=1$. First assume we have $f,g\in L^{q}\left(
B_{1}\right)  $. Define%
\[
u\left(  x\right)  =\int_{B_{1}}P\left(  x,\xi\right)  F\left(  \xi\right)
f\left(  \xi\right)  ^{r}d\xi
\]
for $x\in B_{1}^{+}$. Let $p_{1}$ and $q_{1}$ be the numbers given by%
\[
\frac{n}{p_{1}}=\frac{n-1}{a}+\frac{\left(  n-1\right)  r}{p},\quad\frac
{n}{q_{1}}=\frac{n-1}{a}+\frac{\left(  n-1\right)  r}{q}.
\]
It follows from Proposition \ref{propinterpolation} that%
\begin{align*}
\left\vert u\right\vert _{L^{p_{1}}\left(  B_{1}^{+}\right)  }  & \leq
c\left(  n,p,r,a\right)  \left\vert F\right\vert _{L^{a}\left(  B_{1}\right)
}\left\vert f\right\vert _{L^{p}\left(  B_{1}\right)  }^{r},\\
\left\vert u\right\vert _{L^{q_{1}}\left(  B_{1}^{+}\right)  }  & \leq
c\left(  n,q,r,a\right)  \left\vert F\right\vert _{L^{a}\left(  B_{1}\right)
}\left\vert f\right\vert _{L^{q}\left(  B_{1}\right)  }^{r}.
\end{align*}
Given $0<s<t\leq1/2$. For $\xi\in B_{s}$, a similar calculation as in the
proof of Proposition \ref{prop5.1} shows%
\[
f\left(  \xi\right)  \leq\int_{B_{\frac{s+t}{2}}^{+}}P\left(  x,\xi\right)
U\left(  x\right)  u\left(  x\right)  ^{1/r}dx+\frac{c\left(  n,p,r,a\right)
}{\left(  t-s\right)  ^{n-1}}\left\vert f\right\vert _{L^{p}\left(
B_{1}\right)  }+g\left(  \xi\right)  .
\]
Hence%
\[
\left\vert f\right\vert _{L^{q}\left(  B_{s}\right)  }\leq c\left(
n,q,r,b\right)  \left\vert U\right\vert _{L^{b}\left(  B_{1}^{+}\right)
}\left\vert u\right\vert _{L^{q_{1}}\left(  B_{\frac{s+t}{2}}^{+}\right)
}^{1/r}+\frac{c\left(  n,p,q,r,a\right)  }{\left(  t-s\right)  ^{n-1}%
}\left\vert f\right\vert _{L^{p}\left(  B_{1}\right)  }+\left\vert
g\right\vert _{L^{q}\left(  B_{1/2}\right)  }.
\]
On the other hand, for $x\in B_{\frac{s+t}{2}}^{+}$, we have%
\[
u\left(  x\right)  \leq\int_{B_{t}}P\left(  x,\xi\right)  F\left(  \xi\right)
f\left(  \xi\right)  ^{r}d\xi+\frac{c\left(  n,p,r,a\right)  }{\left(
t-s\right)  ^{n-1}}\left\vert F\right\vert _{L^{a}\left(  B_{1}\right)
}\left\vert f\right\vert _{L^{p}\left(  B_{1}\right)  }^{r}.
\]
Hence%
\[
\left\vert u\right\vert _{L^{q_{1}}\left(  B_{\frac{s+t}{2}}^{+}\right)  }\leq
c\left(  n,q,r,a\right)  \left\vert F\right\vert _{L^{a}\left(  B_{1}\right)
}\left\vert f\right\vert _{L^{q}\left(  B_{t}\right)  }^{r}+\frac{c\left(
n,p,q,r,a\right)  }{\left(  t-s\right)  ^{n-1}}\left\vert F\right\vert
_{L^{a}\left(  B_{1}\right)  }\left\vert f\right\vert _{L^{p}\left(
B_{1}\right)  }^{r}.
\]
Combining the two inequalities together, we see%
\[
\left\vert f\right\vert _{L^{q}\left(  B_{s}\right)  }\leq\frac{1}%
{2}\left\vert f\right\vert _{L^{q}\left(  B_{t}\right)  }+\frac{c\left(
n,p,q,r,a,b\right)  }{\left(  t-s\right)  ^{n-1}}\left\vert f\right\vert
_{L^{p}\left(  B_{1}\right)  }+\left\vert g\right\vert _{L^{q}\left(
B_{1/2}\right)  }%
\]
when $\varepsilon$ is small enough. This implies%
\[
\left\vert f\right\vert _{L^{q}\left(  B_{1/4}\right)  }\leq c\left(
n,p,q,r,a,b\right)  \left(  \left\vert f\right\vert _{L^{p}\left(
B_{1}\right)  }+\left\vert g\right\vert _{L^{q}\left(  B_{1/2}\right)
}\right)  .
\]
With this apriori estimate at hands, we may proceed in the same way as the
proof of Proposition \ref{prop5.1} to get the full conclusion.
\end{proof}

Now we are ready to derive the main results of this section.

\begin{proof}
[Proof of Theorem \ref{thm5.1}]Let $p_{0}=\frac{1}{p-1}$, $f_{0}\left(
\xi\right)  =f\left(  \xi\right)  ^{p-1}$, $u_{0}\left(  x\right)  =\left(
Pf\right)  \left(  x\right)  $, then $0<p_{0}<\infty$, $f_{0}\in
L_{loc}^{p_{0}+1}\left(  \mathbb{R}^{n-1}\right)  $ and%
\[
u_{0}\left(  x\right)  =\int_{\mathbb{R}^{n-1}}P\left(  x,\xi\right)
f_{0}\left(  \xi\right)  ^{p_{0}}d\xi,\quad f_{0}\left(  \xi\right)
=\int_{\mathbb{R}_{+}^{n}}P\left(  x,\xi\right)  u_{0}\left(  x\right)
^{\frac{p_{0}+n}{\left(  n-1\right)  p_{0}}}dx.
\]
For $R>0$, we write $B_{R}=B_{R}^{n-1}$, $B_{R}^{+}=B_{R}^{n}\cap
\mathbb{R}_{+}^{n}$ and%
\begin{align*}
u_{R}\left(  x\right)   & =\int_{\mathbb{R}^{n-1}\backslash B_{R}}P\left(
x,\xi\right)  f_{0}\left(  \xi\right)  ^{p_{0}}d\xi,\\
f_{R}\left(  \xi\right)   & =\int_{\mathbb{R}_{+}^{n}\backslash B_{R}^{+}%
}P\left(  x,\xi\right)  u_{0}\left(  x\right)  ^{\frac{p_{0}+n}{\left(
n-1\right)  p_{0}}}dx.
\end{align*}
First we want to show $u_{0}\in L_{loc}^{\frac{n\left(  p_{0}+1\right)
}{\left(  n-1\right)  p_{0}}}\left(  \overline{\mathbb{R}}_{+}^{n}\right)  $
and $u_{R}\in L^{\frac{n\left(  p_{0}+1\right)  }{\left(  n-1\right)  p_{0}}%
}\left(  B_{R}^{+}\right)  \cap L_{loc}^{\infty}\left(  B_{R}^{+}\cup
B_{R}^{n-1}\right)  $. Indeed, since $f_{0}\in L_{loc}^{p_{0}+1}\left(
\mathbb{R}^{n-1}\right)  $, we see $f_{0}<\infty$ a.e. on $\mathbb{R}^{n-1}$.
This implies $u_{0}<\infty$ a.e. on $\mathbb{R}_{+}^{n}$. Hence there exists
an $x_{0}\in B_{R}^{+}$ such that $u_{0}\left(  x_{0}\right)  <\infty$. It
follows that $\int_{\mathbb{R}^{n-1}\backslash B_{R}}\frac{f_{0}\left(
\xi\right)  ^{p_{0}}}{\left(  \left\vert x_{0}^{\prime}-\xi\right\vert
^{2}+x_{0,n}^{2}\right)  ^{n/2}}d\xi<\infty$ and $\int_{\mathbb{R}%
^{n-1}\backslash B_{R}}\frac{f_{0}\left(  \xi\right)  ^{p_{0}}}{\left\vert
\xi\right\vert ^{n}}d\xi<\infty$. For $0<\theta<1$, $x\in B_{\theta R}^{+}$,
we have%
\[
u_{R}\left(  x\right)  =\int_{\mathbb{R}^{n-1}\backslash B_{R}}P\left(
x,\xi\right)  f_{0}\left(  \xi\right)  ^{p_{0}}d\xi\leq\frac{c\left(
n\right)  R}{\left(  1-\theta\right)  ^{n}}\int_{\mathbb{R}^{n-1}\backslash
B_{R}}\frac{f_{0}\left(  \xi\right)  ^{p_{0}}}{\left\vert \xi\right\vert ^{n}%
}d\xi.
\]
It follows that $u_{R}\in L_{loc}^{\infty}\left(  B_{R}^{+}\cup B_{R}%
^{n-1}\right)  $. Since $\int_{B_{R}}P\left(  \cdot,\xi\right)  f_{0}\left(
\xi\right)  ^{p_{0}}d\xi\in L^{\frac{n\left(  p_{0}+1\right)  }{\left(
n-1\right)  p_{0}}}\left(  \mathbb{R}_{+}^{n}\right)  $, we know $u_{0}\in
L_{loc}^{\frac{n\left(  p_{0}+1\right)  }{\left(  n-1\right)  p_{0}}}\left(
B_{R}^{+}\cup B_{R}^{n-1}\right)  $. By choosing $R$ arbitrarily large, we
deduce that $u_{0}\in L_{loc}^{\frac{n\left(  p_{0}+1\right)  }{\left(
n-1\right)  p_{0}}}\left(  \overline{\mathbb{R}}_{+}^{n}\right)  $ and hence
$u_{R}\in L^{\frac{n\left(  p_{0}+1\right)  }{\left(  n-1\right)  p_{0}}%
}\left(  B_{R}^{+}\right)  $.

Next we want to show $f_{R}\in L^{p_{0}+1}\left(  B_{R}\right)  \cap
L_{loc}^{\infty}\left(  B_{R}\right)  $. Indeed we may find $\xi_{0}\in B_{R}
$ such that $\int_{\mathbb{R}_{+}^{n}}P\left(  x,\xi_{0}\right)  u_{0}\left(
x\right)  ^{\frac{p_{0}+n}{\left(  n-1\right)  p_{0}}}dx<\infty$. This
implies
\[
\int_{\mathbb{R}_{+}^{n}\backslash B_{R}^{+}}\frac{x_{n}}{\left(  \left\vert
x^{\prime}-\xi_{0}\right\vert ^{2}+x_{n}^{2}\right)  ^{n/2}}u_{0}\left(
x\right)  ^{\frac{p_{0}+n}{\left(  n-1\right)  p_{0}}}dx<\infty
\]
and hence $\int_{\mathbb{R}_{+}^{n}\backslash B_{R}^{+}}\frac{x_{n}%
}{\left\vert x\right\vert ^{n}}u_{0}\left(  x\right)  ^{\frac{p_{0}+n}{\left(
n-1\right)  p_{0}}}dx<\infty$. For $0<\theta<1$, $\xi\in B_{\theta R}$, we
have%
\[
f_{R}\left(  \xi\right)  =\int_{\mathbb{R}_{+}^{n}\backslash B_{R}^{+}%
}P\left(  x,\xi\right)  u_{0}\left(  x\right)  ^{\frac{p_{0}+n}{\left(
n-1\right)  p_{0}}}dx\leq\frac{c\left(  n\right)  }{\left(  1-\theta\right)
^{n}}\int_{\mathbb{R}_{+}^{n}\backslash B_{R}^{+}}\frac{x_{n}}{\left\vert
x\right\vert ^{n}}u_{0}\left(  x\right)  ^{\frac{p_{0}+n}{\left(  n-1\right)
p_{0}}}dx
\]
and hence $f_{R}\in L_{loc}^{\infty}\left(  B_{R}\right)  $. To prove the
regularity of $f$, we discuss two cases.

\begin{case}
\label{case5.1}$0<p_{0}\leq\frac{n}{n-1}$.
\end{case}

In this case, we have $\frac{p_{0}+n}{\left(  n-1\right)  p_{0}}>1$. Fix a
number $r$ such that%
\[
1\leq r<\frac{p_{0}+n}{\left(  n-1\right)  p_{0}}\text{ and }r>\frac{1}{p_{0}%
},
\]
then%
\[
f_{0}\left(  \xi\right)  ^{1/r}\leq\left(  \int_{B_{R}^{+}}P\left(
x,\xi\right)  u_{0}\left(  x\right)  ^{\frac{p_{0}+n}{\left(  n-1\right)
p_{0}}}dx\right)  ^{1/r}+f_{R}\left(  \xi\right)  ^{1/r}.
\]
Hence%
\begin{align*}
& u_{0}\left(  x\right) \\
& =\int_{B_{R}}P\left(  x,\xi\right)  f_{0}\left(  \xi\right)  ^{p_{0}-r^{-1}%
}f_{0}\left(  \xi\right)  ^{1/r}d\xi+u_{R}\left(  x\right) \\
& \leq\int_{B_{R}}P\left(  x,\xi\right)  f_{0}\left(  \xi\right)
^{p_{0}-r^{-1}}\left(  \int_{B_{R}^{+}}P\left(  y,\xi\right)  u_{0}\left(
y\right)  ^{\frac{p_{0}+n}{\left(  n-1\right)  p_{0}}-r}u_{0}\left(  y\right)
^{r}dy\right)  ^{1/r}d\xi+v_{R}\left(  x\right)  ,
\end{align*}
here%
\[
v_{R}\left(  x\right)  =\int_{B_{R}}P\left(  x,\xi\right)  f_{0}\left(
\xi\right)  ^{p_{0}-r^{-1}}f_{R}\left(  \xi\right)  ^{1/r}d\xi+u_{R}\left(
x\right)  .
\]
Since $f_{R}\in L^{p_{0}+1}\left(  B_{R}\right)  $, we see $v_{R}\in
L^{\frac{n\left(  p_{0}+1\right)  }{\left(  n-1\right)  p_{0}}}\left(
B_{R}^{+}\right)  $. On the other hand, for $0<\theta<1$, $x\in B_{\theta
R}^{+}$, we have%
\begin{align*}
& \int_{B_{R}}P\left(  x,\xi\right)  f_{0}\left(  \xi\right)  ^{p_{0}-r^{-1}%
}f_{R}\left(  \xi\right)  ^{1/r}d\xi\\
& \leq\left\vert f_{R}\right\vert _{L^{\infty}\left(  B_{\frac{1+\theta}{2}%
R}\right)  }^{1/r}\int_{B_{\frac{1+\theta}{2}R}}P\left(  x,\xi\right)
f_{0}\left(  \xi\right)  ^{p_{0}-r^{-1}}d\xi\\
& +\frac{c\left(  n\right)  }{\left(  1-\theta\right)  ^{n}R^{n-1}}\int
_{B_{R}\backslash B_{\frac{1+\theta}{2}R}}f_{0}\left(  \xi\right)
^{p_{0}-r^{-1}}f_{R}\left(  \xi\right)  ^{1/r}d\xi\\
& \leq\left\vert f_{R}\right\vert _{L^{\infty}\left(  B_{\frac{1+\theta}{2}%
R}\right)  }^{1/r}\int_{B_{\frac{1+\theta}{2}R}}P\left(  x,\xi\right)
f_{0}\left(  \xi\right)  ^{p_{0}-r^{-1}}d\xi+\frac{c\left(  n,p_{0}\right)
}{\left(  1-\theta\right)  ^{n}R^{\frac{\left(  n-1\right)  p_{0}}{p_{0}+1}}%
}\left\vert f_{0}\right\vert _{L^{p_{0}+1}\left(  B_{R}\right)  }^{p_{0}},
\end{align*}
hence $v_{R}\in L_{loc}^{\frac{n\left(  p_{0}+1\right)  }{\left(  n-1\right)
\left(  p_{0}-r^{-1}\right)  }}\left(  B_{R}^{+}\cup B_{R}^{n-1}\right)  $.
Let%
\[
a=\frac{n\left(  p_{0}+1\right)  }{p_{0}+n-\left(  n-1\right)  p_{0}r},\quad
b=\frac{\left(  p_{0}+1\right)  r}{p_{0}r-1}.
\]
Then $\frac{n}{ra}+\frac{n-1}{b}=\frac{1}{r}$ and%
\[
\frac{r}{\frac{n\left(  p_{0}+1\right)  }{\left(  n-1\right)  p_{0}}}+\frac
{1}{a}=\frac{p_{0}+n}{n\left(  p_{0}+1\right)  }<1.
\]
For $\frac{n\left(  p_{0}+1\right)  }{\left(  n-1\right)  p_{0}}%
<q<\frac{n\left(  p_{0}+1\right)  }{\left(  n-1\right)  \left(  p_{0}%
-r^{-1}\right)  }$, we have $\frac{r}{q}+\frac{1}{a}>\frac{1}{n}$. It follows
from Proposition \ref{prop5.1} that $\left.  u_{0}\right\vert _{B_{R/4}^{+}%
}\in L^{q}\left(  B_{R/4}^{+}\right)  $. This implies%
\[
f_{0}\left(  \xi\right)  =\int_{B_{R/4}^{+}}P\left(  x,\xi\right)
u_{0}\left(  x\right)  ^{\frac{p_{0}+n}{\left(  n-1\right)  p_{0}}}%
dx+f_{R/4}\left(  \xi\right)  \leq c\left(  n,q\right)  \left\vert
u_{0}\right\vert _{L^{q}\left(  B_{R/4}^{+}\right)  }^{\frac{p_{0}+n}{\left(
n-1\right)  p_{0}}}+f_{R/4}\left(  \xi\right)
\]
when $q>\frac{n\left(  p_{0}+n\right)  }{\left(  n-1\right)  p_{0}}$. Such a
choice of $q$ is possible since $\frac{n\left(  p_{0}+1\right)  }{\left(
n-1\right)  \left(  p_{0}-r^{-1}\right)  }>\frac{n\left(  p_{0}+n\right)
}{\left(  n-1\right)  p_{0}}$. In particular, we see $\left.  f_{0}\right\vert
_{B_{R/8}}\in L^{\infty}\left(  B_{R/8}\right)  $. Since every point may be
viewed as a center, we see $f_{0}\in L_{loc}^{\infty}\left(  \mathbb{R}%
^{n-1}\right)  $ and hence $u_{0}\in L_{loc}^{\infty}\left(  \overline
{\mathbb{R}}_{+}^{n}\right)  $. For any $R>0$, since%
\[
\int_{\mathbb{R}^{n-1}\backslash B_{R}}\frac{f_{0}\left(  \xi\right)  ^{p_{0}%
}}{\left\vert \xi\right\vert ^{n}}d\xi<\infty\text{ and }\int_{\mathbb{R}%
_{+}^{n}\backslash B_{R}^{+}}\frac{x_{n}}{\left\vert x\right\vert ^{n}}%
u_{0}\left(  x\right)  ^{\frac{p_{0}+n}{\left(  n-1\right)  p_{0}}}dx<\infty,
\]
we see $u_{R}\in C^{\infty}\left(  B_{R}^{+}\cup B_{R}^{n-1}\right)  $ and
$f_{R}\in C^{\infty}\left(  B_{R}\right)  $. It follows that $f_{0}\in
C_{loc}^{\alpha}\left(  \mathbb{R}^{n-1}\right)  $ for $0<\alpha<1$. In
particular, $f_{0}\left(  \xi\right)  >0$ for any $\xi\in\mathbb{R}^{n-1}$.
This implies $u_{0}\in C_{loc}^{\alpha}\left(  \overline{\mathbb{R}}_{+}%
^{n}\right)  $ for any $0<\alpha<1$. Using the fact $\partial_{2}%
\log\left\vert x\right\vert =x_{2}\left\vert x\right\vert ^{-2}$ when $n=2$,
$\partial_{n}\left\vert x\right\vert ^{2-n}=\left(  2-n\right)  x_{n}%
\left\vert x\right\vert ^{-n}$ when $n\geq3$ and the standard potential theory
in \cite[chapter 4]{GT}, it follows from bootstrap method that both $u_{0}$
and $f_{0}$ are smooth. If $f\in L^{p}\left(  \mathbb{R}^{n-1}\right)  $, then
$f_{0}\in L^{p_{0}+1}\left(  \mathbb{R}^{n-1}\right)  $ and $u_{0}\in
L^{\frac{n\left(  p_{0}+1\right)  }{\left(  n-1\right)  p_{0}}}\left(
\mathbb{R}_{+}^{n}\right)  $. If we go back to the proof with this fact and
apply Holder inequality when necessary, we will get $f_{0}\in L^{\infty
}\left(  \mathbb{R}^{n-1}\right)  $ and $u_{0}\in L^{\infty}\left(
\mathbb{R}_{+}^{n}\right)  $. This implies $u_{0}^{\frac{p_{0}+n}{\left(
n-1\right)  p_{0}}}\in L^{s}\left(  \mathbb{R}_{+}^{n}\right)  $ for
$\frac{n\left(  p_{0}+1\right)  }{p_{0}+n}\leq s\leq\infty$. Denote%
\begin{align*}
U  & =\frac{x_{n}}{\left\vert x\right\vert ^{n}}\ast\left(  u_{0}^{\frac
{p_{0}+n}{\left(  n-1\right)  p_{0}}}\chi_{\mathbb{R}_{+}^{n}}\right) \\
& =\left(  \frac{x_{n}}{\left\vert x\right\vert ^{n}}\chi_{B_{1}^{n}}\right)
\ast\left(  u_{0}^{\frac{p_{0}+n}{\left(  n-1\right)  p_{0}}}\chi
_{\mathbb{R}_{+}^{n}}\right)  +\left(  \frac{x_{n}}{\left\vert x\right\vert
^{n}}\chi_{\mathbb{R}^{n}\backslash B_{1}^{n}}\right)  \ast\left(
u_{0}^{\frac{p_{0}+n}{\left(  n-1\right)  p_{0}}}\chi_{\mathbb{R}_{+}^{n}%
}\right)  .
\end{align*}
Since $\frac{x_{n}}{\left\vert x\right\vert ^{n}}\chi_{B_{1}^{n}}\in
L^{\frac{n}{n-1}-\varepsilon}\left(  \mathbb{R}^{n}\right)  $, $\frac{x_{n}%
}{\left\vert x\right\vert ^{n}}\chi_{\mathbb{R}^{n}\backslash B_{1}^{n}}\in
L^{\frac{n}{n-1}+\varepsilon}\left(  \mathbb{R}^{n}\right)  $ and
$\frac{n\left(  p_{0}+1\right)  }{p_{0}+n}<n$, we see $U$ is continuous and
$U\left(  x\right)  \rightarrow0$ as $\left\vert x\right\vert \rightarrow
\infty$. It follows from $f_{0}=c\left(  n\right)  \left.  U\right\vert
_{\mathbb{R}^{n-1}}$ that $f_{0}\left(  \xi\right)  \rightarrow0$ as
$\left\vert \xi\right\vert \rightarrow\infty$.

\begin{case}
\label{case5.2}$\frac{n}{n-1}\leq p_{0}<\infty$.
\end{case}

In this case, we fix a number $r$ such that%
\[
1\leq r\leq p_{0}\text{ and }r\geq\frac{\left(  n-1\right)  p_{0}}{p_{0}+n},
\]
then%
\[
u_{0}\left(  x\right)  ^{1/r}\leq\left(  \int_{B_{R}}P\left(  x,\xi\right)
f_{0}\left(  \xi\right)  ^{p_{0}}d\xi\right)  ^{1/r}+u_{R}\left(  x\right)
^{1/r}.
\]
Hence%
\[
f_{0}\left(  \xi\right)  \leq\int_{B_{R}^{+}}P\left(  x,\xi\right)
u_{0}\left(  x\right)  ^{\frac{p_{0}+n}{\left(  n-1\right)  p_{0}}-r^{-1}%
}\left(  \int_{B_{R}}P\left(  x,\zeta\right)  f_{0}\left(  \zeta\right)
^{p_{0}-r}f_{0}\left(  \zeta\right)  ^{r}d\zeta\right)  ^{1/r}dx+g_{R}\left(
\xi\right)  ,
\]
here%
\[
g_{R}\left(  \xi\right)  =\int_{B_{R}^{+}}P\left(  x,\xi\right)  u_{0}\left(
x\right)  ^{\frac{p_{0}+n}{\left(  n-1\right)  p_{0}}-r^{-1}}u_{R}\left(
x\right)  ^{1/r}dx+f_{R}\left(  \xi\right)  .
\]
Since $u_{R}\in L^{\frac{n\left(  p_{0}+1\right)  }{\left(  n-1\right)  p_{0}%
}}\left(  B_{R}^{+}\right)  $, we see $g_{R}\in L^{p_{0}+1}\left(
B_{R}\right)  $. On the other hand, for $0<\theta<1$, $\xi\in B_{\theta R}$,
we have%
\begin{align*}
& \int_{B_{R}^{+}}P\left(  x,\xi\right)  u_{0}\left(  x\right)  ^{\frac
{p_{0}+n}{\left(  n-1\right)  p_{0}}-r^{-1}}u_{R}\left(  x\right)  ^{1/r}dx\\
& \leq\left\vert u_{R}\right\vert _{L^{\infty}\left(  B_{\frac{1+\theta}{2}%
R}^{+}\right)  }^{1/r}\int_{B_{\frac{1+\theta}{2}R}^{+}}P\left(  x,\xi\right)
u_{0}\left(  x\right)  ^{\frac{p_{0}+n}{\left(  n-1\right)  p_{0}}-r^{-1}}dx\\
& +\frac{c\left(  n\right)  }{\left(  1-\theta\right)  ^{n}R^{n-1}}\int
_{B_{R}^{+}\backslash B_{\frac{1+\theta}{2}R}^{+}}u_{0}\left(  x\right)
^{\frac{p_{0}+n}{\left(  n-1\right)  p_{0}}-r^{-1}}u_{R}\left(  x\right)
^{1/r}dx\\
& \leq\left\vert u_{R}\right\vert _{L^{\infty}\left(  B_{\frac{1+\theta}{2}%
R}^{+}\right)  }^{1/r}\int_{B_{\frac{1+\theta}{2}R}^{+}}P\left(  x,\xi\right)
u_{0}\left(  x\right)  ^{\frac{p_{0}+n}{\left(  n-1\right)  p_{0}}-r^{-1}}dx\\
& +\frac{c\left(  n,p_{0}\right)  }{\left(  1-\theta\right)  ^{n}R^{\frac
{n-1}{p_{0}+1}}}\left\vert u_{0}\right\vert _{L^{\frac{n\left(  p_{0}%
+1\right)  }{\left(  n-1\right)  p_{0}}}\left(  B_{R}^{+}\right)  }%
^{\frac{p_{0}+n}{\left(  n-1\right)  p_{0}}},
\end{align*}
hence $g_{R}\in L_{loc}^{q}\left(  B_{R}\right)  $ for any $q<\infty$. Let%
\[
a=\frac{p_{0}+1}{p_{0}-r},\quad b=\frac{n\left(  p_{0}+1\right)  r}{\left(
p_{0}+n\right)  r-\left(  n-1\right)  p_{0}},
\]
then $\frac{n-1}{ra}+\frac{n}{b}=1$, $\frac{r}{p_{0}+1}+\frac{1}{a}%
=\frac{p_{0}}{p_{0}+1}\in\left(  0,1\right)  $. For any $p_{0}+1<q<\infty$, it
follows from Proposition \ref{prop5.2} that when $R$ is small enough, we have
$f_{0}\in L^{q}\left(  B_{R/4}\right)  $. Since every point can be viewed as a
center, we see $f_{0}\in L_{loc}^{q}\left(  \mathbb{R}^{n-1}\right)  $ and
hence $u_{0}\in L_{loc}^{\frac{nq}{n-1}}\left(  \overline{\mathbb{R}}_{+}%
^{n}\right)  $. Using the equations of $f_{0}$ and $u_{0}$, we see $f_{0}\in
L_{loc}^{\infty}\left(  \mathbb{R}^{n-1}\right)  $ and $u_{0}\in
L_{loc}^{\infty}\left(  \overline{\mathbb{R}}_{+}^{n}\right)  $. Now the
arguments in Case \ref{case5.1} tell us $f_{0}\in C^{\infty}\left(
\mathbb{R}^{n-1}\right)  $ and $u_{0}\in C^{\infty}\left(  \overline
{\mathbb{R}}_{+}^{n}\right)  $, moreover, $f_{0}\left(  \xi\right)
\rightarrow0$ as $\left\vert \xi\right\vert \rightarrow\infty$ under the
assumption $f\in L^{p}\left(  \mathbb{R}^{n-1}\right)  $.
\end{proof}

\section{Radial symmetry of nonnegative critical functions\label{sec6}}

In this section we will study the symmetry property of the nonnegative
critical functions of the variational problem (\ref{eq1.8}). We will show any
nonnegative critical functions are radial symmetric with respect to some
point. As explained at the beginning of Section \ref{sec5}, (\ref{eq1.9}) may
be viewed as an integral system which is very similar to the integral systems
related to the Hardy-Littlewood-Sobolev inequalities. For the latter one, the
radial symmetry of nonnegative solution for some special cases were solved in
\cite{CLO1,CLO2,L}. In particular, in \cite{CLO2} an integral version of the
method of moving planes (\cite{GNN}) was introduced and later applied in
\cite{CLO1} to resolve the symmetry problems for some cases of the integral
systems related to Hardy-Littlewood-Sobolev inequalities. In \cite{Hn}, some
new observations were added and all the cases for the symmetry of the
solutions to the systems were resolved. We will apply these new observations
to (\ref{eq1.9}).

\begin{theorem}
\label{thm6.1}Assume $1<p<\infty$, $n\geq2$, $f\in L^{p}\left(  \mathbb{R}%
^{n-1}\right)  $ is nonnegative, not identically zero and it satisfies%
\[
f\left(  \xi\right)  ^{p-1}=\int_{\mathbb{R}_{+}^{n}}P\left(  x,\xi\right)
\left(  Pf\right)  \left(  x\right)  ^{\frac{np}{n-1}-1}dx,
\]
then $f\in C^{\infty}\left(  \mathbb{R}^{n-1}\right)  $, moreover $f$ is
radial symmetric with respect to some point and strictly decreasing along the
radial direction.
\end{theorem}

For the case $n\geq3$, $p=\frac{2\left(  n-1\right)  }{n-2}$, the
Euler-Lagrange equation has conformal invariance property and we may weaken
the assumption a little bit.

\begin{proposition}
\label{prop6.1}Assume $n\geq3$, $f\in L_{loc}^{\frac{2\left(  n-1\right)
}{n-2}}\left(  \mathbb{R}^{n-1}\right)  $ is nonnegative, not identically zero
and it satisfies%
\[
f\left(  \xi\right)  ^{\frac{n}{n-2}}=\int_{\mathbb{R}_{+}^{n}}P\left(
x,\xi\right)  \left(  Pf\right)  \left(  x\right)  ^{\frac{n+2}{n-2}}dx,
\]
then for some $\lambda>0$ and $\xi_{0}\in\mathbb{R}^{n-1}$, we have%
\[
f\left(  \xi\right)  =c\left(  n\right)  \left(  \frac{\lambda}{\lambda
^{2}+\left\vert \xi-\xi_{0}\right\vert ^{2}}\right)  ^{\frac{n-2}{2}}.
\]

\end{proposition}

In the proof of these symmetry results, we will need the following basic
inequality: assume $0<\theta\leq1$, $a\geq b\geq0$, $c\geq0$, then%
\[
\left(  a+c\right)  ^{\theta}-\left(  b+c\right)  ^{\theta}\leq a^{\theta
}-b^{\theta}.
\]

For $\sigma\in\mathbb{R}^{m}$ and $s>0$, we denote%
\[
\left\vert \sigma\right\vert _{l^{s}}=\left(  \sum_{i=1}^{m}\left\vert
\sigma_{i}\right\vert ^{s}\right)  ^{1/s}.
\]

\begin{proof}
[Proof of Theorem \ref{thm6.1}]By Theorem \ref{thm5.1} we know $f\in
C^{\infty}\left(  \mathbb{R}^{n-1}\right)  $ and $f\left(  \xi\right)
\rightarrow0$ as $\left\vert \xi\right\vert \rightarrow\infty$. Let
$q=\frac{1}{p-1}$, $g\left(  \xi\right)  =f\left(  \xi\right)  ^{p-1}$,
$v\left(  x\right)  =\left(  Pf\right)  \left(  x\right)  $, then $0<q<\infty
$, $g\in L^{q+1}\left(  \mathbb{R}^{n-1}\right)  $, $v\in L^{\frac{n\left(
q+1\right)  }{\left(  n-1\right)  q}}\left(  \mathbb{R}_{+}^{n}\right)  $ and%
\[
v\left(  x\right)  =\int_{\mathbb{R}^{n-1}}P\left(  x,\xi\right)  g\left(
\xi\right)  ^{q}d\xi,\quad g\left(  \xi\right)  =\int_{\mathbb{R}_{+}^{n}%
}P\left(  x,\xi\right)  v\left(  x\right)  ^{\frac{q+n}{\left(  n-1\right)
q}}dx.
\]
For $\lambda\in\mathbb{R}$, denote%
\[
H_{\lambda}=\left\{  \xi\in\mathbb{R}^{n-1}:\xi_{1}<\lambda\right\}  ,\quad
Q_{\lambda}=\left\{  x\in\mathbb{R}_{+}^{n}:x_{1}<\lambda\right\}  .
\]
For $\xi\in\mathbb{R}^{n-1}$, $\xi=\left(  \xi_{1},\xi^{\prime\prime}\right)
$, denote $\xi_{\lambda}=\left(  2\lambda-\xi_{1},\xi^{\prime\prime}\right)
$. For $x\in\mathbb{R}^{n}$, $x=\left(  x_{1},x^{\prime\prime}\right)  $,
denote $x_{\lambda}=\left(  2\lambda-x_{1},x^{\prime\prime}\right)  $. Define
$g_{\lambda}\left(  \xi\right)  =g\left(  \xi_{\lambda}\right)  $,
$v_{\lambda}\left(  x\right)  =v\left(  x_{\lambda}\right)  $ and%
\[
\mathcal{B}_{\lambda}^{g}=\left\{  \xi\in H_{\lambda}:g_{\lambda}\left(
\xi\right)  >g\left(  \xi\right)  \right\}  ,\quad\mathcal{B}_{\lambda}%
^{v}=\left\{  x\in Q_{\lambda}:v_{\lambda}\left(  x\right)  >v\left(
x\right)  \right\}  .
\]
By a simple change of variable, we see%
\begin{align*}
v\left(  x\right)   & =\int_{H_{\lambda}}P\left(  x,\xi\right)  g\left(
\xi\right)  ^{q}d\xi+\int_{H_{\lambda}}P\left(  x_{\lambda},\xi\right)
g\left(  \xi_{\lambda}\right)  ^{q}d\xi,\\
g\left(  \xi\right)   & =\int_{Q_{\lambda}}P\left(  x,\xi\right)  v\left(
x\right)  ^{\frac{q+n}{\left(  n-1\right)  q}}dx+\int_{Q_{\lambda}}P\left(
x_{\lambda},\xi\right)  v\left(  x_{\lambda}\right)  ^{\frac{q+n}{\left(
n-1\right)  q}}dx.
\end{align*}

\begin{case}
\label{case6.1}$0<q\leq\frac{n}{n-1}$.
\end{case}

In this case, we choose a number $r$ such that%
\[
1\leq r\leq\frac{q+n}{\left(  n-1\right)  q}\text{ and }q^{-1}<r.
\]
We have%
\[
v\left(  x_{\lambda}\right)  -v\left(  x\right)  =\int_{H_{\lambda}}\left(
P\left(  x,\xi\right)  -P\left(  x_{\lambda},\xi\right)  \right)  \left(
g\left(  \xi_{\lambda}\right)  ^{q}-g\left(  \xi\right)  ^{q}\right)  d\xi.
\]
Hence for $x\in\mathcal{B}_{\lambda}^{v}$,%
\begin{align*}
0  & \leq v\left(  x_{\lambda}\right)  -v\left(  x\right) \\
& \leq\int_{\mathcal{B}_{\lambda}^{g}}\left(  P\left(  x,\xi\right)  -P\left(
x_{\lambda},\xi\right)  \right)  \left(  g\left(  \xi_{\lambda}\right)
^{q}-g\left(  \xi\right)  ^{q}\right)  d\xi\\
& \leq qr\int_{\mathcal{B}_{\lambda}^{g}}P\left(  x,\xi\right)  g\left(
\xi_{\lambda}\right)  ^{q-r^{-1}}\left(  g\left(  \xi_{\lambda}\right)
^{1/r}-g\left(  \xi\right)  ^{1/r}\right)  d\xi.
\end{align*}
It follows that%
\begin{align*}
\left\vert v_{\lambda}-v\right\vert _{L^{\frac{n\left(  q+1\right)  }{\left(
n-1\right)  q}}\left(  \mathcal{B}_{\lambda}^{v}\right)  }  & \leq c\left(
n,q,r\right)  \left\vert g_{\lambda}^{q-r^{-1}}\left(  g_{\lambda}%
^{1/r}-g^{1/r}\right)  \right\vert _{L^{\frac{q+1}{q}}\left(  \mathcal{B}%
_{\lambda}^{g}\right)  }\\
& \leq c\left(  n,q,r\right)  \left\vert g_{\lambda}\right\vert _{L^{q+1}%
\left(  \mathcal{B}_{\lambda}^{g}\right)  }^{q-r^{-1}}\left\vert g_{\lambda
}^{1/r}-g^{1/r}\right\vert _{L^{\left(  q+1\right)  r}\left(  \mathcal{B}%
_{\lambda}^{g}\right)  }.
\end{align*}
On the other hand, for $\xi\in\mathcal{B}_{\lambda}^{g}$, we have%
\begin{align*}
g\left(  \xi_{\lambda}\right)   & =\int_{\mathcal{B}_{\lambda}^{v}}P\left(
x,\xi\right)  v\left(  x_{\lambda}\right)  ^{\frac{q+n}{\left(  n-1\right)
q}}dx+\int_{\mathcal{B}_{\lambda}^{v}}P\left(  x_{\lambda},\xi\right)
v\left(  x\right)  ^{\frac{q+n}{\left(  n-1\right)  q}}dx\\
& +\int_{Q_{\lambda}\backslash\mathcal{B}_{\lambda}^{v}}P\left(  x,\xi\right)
v\left(  x_{\lambda}\right)  ^{\frac{q+n}{\left(  n-1\right)  q}}%
dx+\int_{Q_{\lambda}\backslash\mathcal{B}_{\lambda}^{v}}P\left(  x_{\lambda
},\xi\right)  v\left(  x\right)  ^{\frac{q+n}{\left(  n-1\right)  q}}dx\\
& \leq\int_{\mathcal{B}_{\lambda}^{v}}P\left(  x,\xi\right)  v\left(
x_{\lambda}\right)  ^{\frac{q+n}{\left(  n-1\right)  q}}dx+\int_{\mathcal{B}%
_{\lambda}^{v}}P\left(  x_{\lambda},\xi\right)  v\left(  x\right)
^{\frac{q+n}{\left(  n-1\right)  q}}dx\\
& +\int_{Q_{\lambda}\backslash\mathcal{B}_{\lambda}^{v}}P\left(  x,\xi\right)
v\left(  x\right)  ^{\frac{q+n}{\left(  n-1\right)  q}}dx+\int_{Q_{\lambda
}\backslash\mathcal{B}_{\lambda}^{v}}P\left(  x_{\lambda},\xi\right)  v\left(
x_{\lambda}\right)  ^{\frac{q+n}{\left(  n-1\right)  q}}dx.
\end{align*}
Since%
\begin{align*}
g\left(  \xi\right)   & =\int_{\mathcal{B}_{\lambda}^{v}}P\left(
x,\xi\right)  v\left(  x\right)  ^{\frac{q+n}{\left(  n-1\right)  q}}%
dx+\int_{\mathcal{B}_{\lambda}^{v}}P\left(  x_{\lambda},\xi\right)  v\left(
x_{\lambda}\right)  ^{\frac{q+n}{\left(  n-1\right)  q}}dx\\
& +\int_{Q_{\lambda}\backslash\mathcal{B}_{\lambda}^{v}}P\left(  x,\xi\right)
v\left(  x\right)  ^{\frac{q+n}{\left(  n-1\right)  q}}dx+\int_{Q_{\lambda
}\backslash\mathcal{B}_{\lambda}^{v}}P\left(  x_{\lambda},\xi\right)  v\left(
x_{\lambda}\right)  ^{\frac{q+n}{\left(  n-1\right)  q}}dx,
\end{align*}
we see%
\begin{align*}
& g\left(  \xi_{\lambda}\right)  ^{1/r}-g\left(  \xi\right)  ^{1/r}\\
& \leq\left(  \int_{\mathcal{B}_{\lambda}^{v}}P\left(  x,\xi\right)  v\left(
x_{\lambda}\right)  ^{\frac{q+n}{\left(  n-1\right)  q}}dx+\int_{\mathcal{B}%
_{\lambda}^{v}}P\left(  x_{\lambda},\xi\right)  v\left(  x\right)
^{\frac{q+n}{\left(  n-1\right)  q}}dx\right)  ^{1/r}\\
& -\left(  \int_{\mathcal{B}_{\lambda}^{v}}P\left(  x,\xi\right)  v\left(
x\right)  ^{\frac{q+n}{\left(  n-1\right)  q}}dx+\int_{\mathcal{B}_{\lambda
}^{v}}P\left(  x_{\lambda},\xi\right)  v\left(  x_{\lambda}\right)
^{\frac{q+n}{\left(  n-1\right)  q}}dx\right)  ^{1/r}\\
& =\left(  \int_{\mathcal{B}_{\lambda}^{v}}\left\vert \left(  P\left(
x,\xi\right)  ^{1/r}v\left(  x_{\lambda}\right)  ^{\frac{q+n}{\left(
n-1\right)  qr}},P\left(  x_{\lambda},\xi\right)  ^{1/r}v\left(  x\right)
^{\frac{q+n}{\left(  n-1\right)  qr}}\right)  \right\vert _{l^{r}}%
^{r}dx\right)  ^{1/r}\\
& -\left(  \int_{\mathcal{B}_{\lambda}^{v}}\left\vert \left(  P\left(
x,\xi\right)  ^{1/r}v\left(  x\right)  ^{\frac{q+n}{\left(  n-1\right)  qr}%
},P\left(  x_{\lambda},\xi\right)  ^{1/r}v\left(  x_{\lambda}\right)
^{\frac{q+n}{\left(  n-1\right)  qr}}\right)  \right\vert _{l^{r}}%
^{r}dx\right)  ^{1/r}\\
& \leq\left(  \int_{\mathcal{B}_{\lambda}^{v}}\left\vert \left(  I,II\right)
\right\vert _{l^{r}}^{r}dx\right)  ^{1/r}.
\end{align*}
Here%
\begin{align*}
I  & =P\left(  x,\xi\right)  ^{1/r}\left(  v\left(  x_{\lambda}\right)
^{\frac{q+n}{\left(  n-1\right)  qr}}-v\left(  x\right)  ^{\frac{q+n}{\left(
n-1\right)  qr}}\right)  ,\\
II  & =P\left(  x_{\lambda},\xi\right)  ^{1/r}\left(  v\left(  x\right)
^{\frac{q+n}{\left(  n-1\right)  qr}}-v\left(  x_{\lambda}\right)
^{\frac{q+n}{\left(  n-1\right)  qr}}\right)  .
\end{align*}
Hence%
\begin{align*}
& g\left(  \xi_{\lambda}\right)  ^{1/r}-g\left(  \xi\right)  ^{1/r}\\
& \leq2\left(  \int_{\mathcal{B}_{\lambda}^{v}}P\left(  x,\xi\right)  \left(
v\left(  x_{\lambda}\right)  ^{\frac{q+n}{\left(  n-1\right)  qr}}-v\left(
x\right)  ^{\frac{q+n}{\left(  n-1\right)  qr}}\right)  ^{r}dx\right)
^{1/r}\\
& \leq\frac{2\left(  q+n\right)  }{\left(  n-1\right)  qr}\left(
\int_{\mathcal{B}_{\lambda}^{v}}P\left(  x,\xi\right)  v\left(  x_{\lambda
}\right)  ^{\frac{q+n}{\left(  n-1\right)  q}-r}\left(  v\left(  x_{\lambda
}\right)  -v\left(  x\right)  \right)  ^{r}dx\right)  ^{1/r}.
\end{align*}
This implies%
\begin{align*}
& \left\vert g_{\lambda}^{1/r}-g^{1/r}\right\vert _{L^{\left(  q+1\right)
r}\left(  \mathcal{B}_{\lambda}^{g}\right)  }\\
& \leq\frac{2\left(  q+n\right)  }{\left(  n-1\right)  qr}\left\vert
\int_{\mathcal{B}_{\lambda}^{v}}P\left(  x,\xi\right)  v\left(  x_{\lambda
}\right)  ^{\frac{q+n}{\left(  n-1\right)  q}-r}\left(  v\left(  x_{\lambda
}\right)  -v\left(  x\right)  \right)  ^{r}dx\right\vert _{L^{q+1}\left(
\mathcal{B}_{\lambda}^{g}\right)  }^{1/r}\\
& \leq c\left(  n,q,r\right)  \left\vert v_{\lambda}^{\frac{q+n}{\left(
n-1\right)  q}-r}\left(  v_{\lambda}-v\right)  ^{r}\right\vert _{L^{\frac
{n\left(  q+1\right)  }{q+n}}\left(  \mathcal{B}_{\lambda}^{v}\right)  }%
^{1/r}\\
& \leq c\left(  n,q,r\right)  \left\vert v_{\lambda}\right\vert _{L^{\frac
{n\left(  q+1\right)  }{\left(  n-1\right)  q}}\left(  \mathcal{B}_{\lambda
}^{v}\right)  }^{\frac{q+n}{\left(  n-1\right)  qr}-1}\left\vert v_{\lambda
}-v\right\vert _{L^{\frac{n\left(  q+1\right)  }{\left(  n-1\right)  q}%
}\left(  \mathcal{B}_{\lambda}^{v}\right)  }.
\end{align*}
It follows from the two inequalities above that%
\[
\left\vert g_{\lambda}^{1/r}-g^{1/r}\right\vert _{L^{\left(  q+1\right)
r}\left(  \mathcal{B}_{\lambda}^{g}\right)  }\leq c\left(  n,q,r\right)
\left\vert v\right\vert _{L^{\frac{n\left(  q+1\right)  }{\left(  n-1\right)
q}}\left(  \mathbb{R}_{+}^{n}\right)  }^{\frac{q+n}{\left(  n-1\right)  qr}%
-1}\left\vert g\right\vert _{L^{q+1}\left(  2\lambda e_{1}^{\prime
}-\mathcal{B}_{\lambda}^{g}\right)  }^{q-r^{-1}}\left\vert g_{\lambda}%
^{1/r}-g^{1/r}\right\vert _{L^{\left(  q+1\right)  r}\left(  \mathcal{B}%
_{\lambda}^{g}\right)  }.
\]
Here $e_{1}=\left(  1,0,\cdots,0\right)  $. After these preparations, we will
use the method of moving planes to prove the radial symmetry of $g$ and hence
$f$.

First, we have to show it is possible to start. Indeed, for $\lambda$ large
enough, $\left\vert g\right\vert _{L^{q+1}\left(  2\lambda e_{1}^{\prime
}-\mathcal{B}_{\lambda}^{v}\right)  }$ can be arbitrarily small, this implies%
\[
\left\vert g_{\lambda}^{1/r}-g^{1/r}\right\vert _{L^{\left(  q+1\right)
r}\left(  \mathcal{B}_{\lambda}^{g}\right)  }\leq\frac{1}{2}\left\vert
g_{\lambda}^{1/r}-g^{1/r}\right\vert _{L^{\left(  q+1\right)  r}\left(
\mathcal{B}_{\lambda}^{g}\right)  }%
\]
and hence $\left\vert g_{\lambda}^{1/r}-g^{1/r}\right\vert _{L^{\left(
q+1\right)  r}\left(  \mathcal{B}_{\lambda}^{g}\right)  }=0$. It follows that
$\mathcal{B}_{\lambda}^{g}=\emptyset$ when $\lambda$ is large enough.

Next we let $\lambda_{0}=\inf\left\{  \lambda\in\mathbb{R}:\mathcal{B}%
_{\lambda^{\prime}}^{g}=\emptyset\text{ for all }\lambda^{\prime}\geq
\lambda\right\}  $. It follows from the fact $g\left(  \xi\right)
\rightarrow0$ as $\left\vert \xi\right\vert \rightarrow\infty$ and $g\left(
\xi\right)  >0$ for all $\xi\in\mathbb{R}^{n-1}$ that $\lambda_{0} $ must be a
finite number. By the definition of $\lambda_{0}$ we know $g_{\lambda_{0}%
}\left(  \xi\right)  \leq g\left(  \xi\right)  $ for $\xi\in H_{\lambda_{0}}$.
We claim that $g_{\lambda_{0}}=g$. Indeed if this is not the case, then since%
\[
v_{\lambda_{0}}\left(  x\right)  -v\left(  x\right)  =\int_{H_{\lambda_{0}}%
}\left(  P\left(  x,\xi\right)  -P\left(  x_{\lambda_{0}},\xi\right)  \right)
\left(  g_{\lambda_{0}}\left(  \xi\right)  ^{q}-g\left(  \xi\right)
^{q}\right)  d\xi
\]
and%
\[
g_{\lambda_{0}}\left(  \xi\right)  -g\left(  \xi\right)  =\int_{Q_{\lambda
_{0}}}\left(  P\left(  x,\xi\right)  -P\left(  x_{\lambda_{0}},\xi\right)
\right)  \left(  v_{\lambda_{0}}\left(  x\right)  ^{\frac{q+n}{\left(
n-1\right)  q}}-v\left(  x\right)  ^{\frac{q+n}{\left(  n-1\right)  q}%
}\right)  dx,
\]
we get $g_{\lambda_{0}}\left(  \xi\right)  <g\left(  \xi\right)  $ for $\xi\in
H_{\lambda_{0}}$. It follows that $\chi_{2\lambda e_{1}^{\prime}%
-\mathcal{B}_{\lambda}^{g}}\rightarrow0$ a.e. as $\lambda\uparrow\lambda_{0}$.
By dominated convergence theorem we have $\left\vert g\right\vert
_{L^{q+1}\left(  2\lambda e_{1}^{\prime}-\mathcal{B}_{\lambda}^{g}\right)
}\rightarrow0$ as $\lambda\uparrow\lambda_{0}$. It implies%
\[
\left\vert g_{\lambda}^{1/r}-g^{1/r}\right\vert _{L^{\left(  q+1\right)
r}\left(  \mathcal{B}_{\lambda}^{g}\right)  }\leq\frac{1}{2}\left\vert
g_{\lambda}^{1/r}-g^{1/r}\right\vert _{L^{\left(  q+1\right)  r}\left(
\mathcal{B}_{\lambda}^{g}\right)  }%
\]
when $\lambda$ is very close to $\lambda_{0}$ and hence $\mathcal{B}_{\lambda
}^{g}=\emptyset$. This contradicts with the choice of $\lambda_{0}$. Hence
when the moving process stops, we must have symmetry. Moreover we claim that
$g_{\lambda}\left(  \xi\right)  <g\left(  \xi\right)  $ for $\xi\in
H_{\lambda}$ when $\lambda>\lambda_{0}$. Indeed for any $\lambda>\lambda_{0}$
we can not have $g_{\lambda}=g$ because otherwise $g$ is periodic in the first
direction and can not lie in $L^{q+1}\left(  \mathbb{R}^{n-1}\right)  $. Hence
$g_{\lambda}<g$ in $H_{\lambda}$.

By translation, we may assume $g\left(  0\right)  =\max_{\xi\in\mathbb{R}%
^{n-1}}g\left(  \xi\right)  $, then it follows that the moving plane process
from any direction must stop at the origin. Hence $g$ must be radial symmetric
and strictly decreasing in the radial direction.

\begin{case}
\label{case6.2}$\frac{n}{n-1}\leq q<\infty$.
\end{case}

In this case, we choose a number $r$ such that%
\[
1\leq r<q\text{ and }\frac{\left(  n-1\right)  q}{q+n}\leq r.
\]
We have%
\[
g\left(  \xi_{\lambda}\right)  -g\left(  \xi\right)  =\int_{Q_{\lambda}%
}\left(  P\left(  x,\xi\right)  -P\left(  x_{\lambda},\xi\right)  \right)
\left(  v\left(  x_{\lambda}\right)  ^{\frac{q+n}{\left(  n-1\right)  q}%
}-v\left(  x\right)  ^{\frac{q+n}{\left(  n-1\right)  q}}\right)  dx.
\]
Hence for $\xi\in\mathcal{B}_{\lambda}^{g}$,%
\begin{align*}
0  & \leq g\left(  \xi_{\lambda}\right)  -g\left(  \xi\right) \\
& \leq\int_{\mathcal{B}_{\lambda}^{v}}\left(  P\left(  x,\xi\right)  -P\left(
x,\xi_{\lambda}\right)  \right)  \left(  v\left(  x_{\lambda}\right)
^{\frac{q+n}{\left(  n-1\right)  q}}-v\left(  x\right)  ^{\frac{q+n}{\left(
n-1\right)  q}}\right)  dx\\
& \leq\frac{\left(  q+n\right)  r}{\left(  n-1\right)  q}\int_{\mathcal{B}%
_{\lambda}^{v}}P\left(  x,\xi\right)  v\left(  x_{\lambda}\right)
^{\frac{q+n}{\left(  n-1\right)  q}-r^{-1}}\left(  v\left(  x_{\lambda
}\right)  ^{1/r}-v\left(  x\right)  ^{1/r}\right)  dx.
\end{align*}
It follows that%
\begin{align*}
& \left\vert g_{\lambda}-g\right\vert _{L^{q+1}\left(  \mathcal{B}_{\lambda
}^{g}\right)  }\\
& \leq c\left(  n,q,r\right)  \left\vert v_{\lambda}^{\frac{q+n}{\left(
n-1\right)  q}-r^{-1}}\left(  v_{\lambda}^{1/r}-v^{1/r}\right)  \right\vert
_{L^{\frac{n\left(  q+1\right)  }{q+n}}\left(  \mathcal{B}_{\lambda}%
^{v}\right)  }\\
& =c\left(  n,q,r\right)  \left\vert v_{\lambda}\right\vert _{L^{\frac
{n\left(  q+1\right)  }{\left(  n-1\right)  q}}\left(  \mathcal{B}_{\lambda
}^{v}\right)  }^{\frac{q+n}{\left(  n-1\right)  q}-r^{-1}}\left\vert
v_{\lambda}^{1/r}-v^{1/r}\right\vert _{L^{\frac{n\left(  q+1\right)
r}{\left(  n-1\right)  q}}\left(  \mathcal{B}_{\lambda}^{v}\right)  }.
\end{align*}
On the other hand, for $x\in\mathcal{B}_{\lambda}^{v}$, we have%
\begin{align*}
v\left(  x_{\lambda}\right)   & \leq\int_{\mathcal{B}_{\lambda}^{g}}P\left(
x,\xi\right)  g\left(  \xi_{\lambda}\right)  ^{q}d\xi+\int_{\mathcal{B}%
_{\lambda}^{g}}P\left(  x_{\lambda},\xi\right)  g\left(  \xi\right)  ^{q}%
d\xi\\
& +\int_{H_{\lambda}\backslash\mathcal{B}_{\lambda}^{g}}P\left(  x,\xi\right)
g\left(  \xi\right)  ^{q}d\xi+\int_{H_{\lambda}\backslash\mathcal{B}_{\lambda
}^{g}}P\left(  x_{\lambda},\xi\right)  g\left(  \xi_{\lambda}\right)  ^{q}%
d\xi.
\end{align*}
Since%
\begin{align*}
v\left(  x\right)   & =\int_{\mathcal{B}_{\lambda}^{g}}P\left(  x,\xi\right)
g\left(  \xi\right)  ^{q}d\xi+\int_{\mathcal{B}_{\lambda}^{g}}P\left(
x_{\lambda},\xi\right)  g\left(  \xi_{\lambda}\right)  ^{q}d\xi\\
& +\int_{H_{\lambda}\backslash\mathcal{B}_{\lambda}^{g}}P\left(  x,\xi\right)
g\left(  \xi\right)  ^{q}d\xi+\int_{H_{\lambda}\backslash\mathcal{B}_{\lambda
}^{g}}P\left(  x_{\lambda},\xi\right)  g\left(  \xi_{\lambda}\right)  ^{q}%
d\xi,
\end{align*}
we see%
\begin{align*}
0  & \leq v\left(  x_{\lambda}\right)  ^{1/r}-v\left(  x\right)  ^{1/r}\\
& \leq\left(  \int_{\mathcal{B}_{\lambda}^{g}}P\left(  x,\xi\right)  g\left(
\xi_{\lambda}\right)  ^{q}d\xi+\int_{\mathcal{B}_{\lambda}^{g}}P\left(
x_{\lambda},\xi\right)  g\left(  \xi\right)  ^{q}d\xi\right)  ^{1/r}\\
& -\left(  \int_{\mathcal{B}_{\lambda}^{g}}P\left(  x,\xi\right)  g\left(
\xi\right)  ^{q}d\xi+\int_{\mathcal{B}_{\lambda}^{g}}P\left(  x_{\lambda}%
,\xi\right)  g\left(  \xi_{\lambda}\right)  ^{q}d\xi\right)  ^{1/r}\\
& \leq\left(  \int_{\mathcal{B}_{\lambda}^{g}}\left\vert \left(  P\left(
x,\xi\right)  ^{1/r}\left(  g\left(  \xi_{\lambda}\right)  ^{q/r}-g\left(
\xi\right)  ^{q/r}\right)  ,P\left(  x_{\lambda},\xi\right)  ^{1/r}\left(
g\left(  \xi\right)  ^{q/r}-g\left(  \xi_{\lambda}\right)  ^{q/r}\right)
\right)  \right\vert _{l^{r}}^{r}d\xi\right)  ^{1/r}\\
& \leq2\left(  \int_{\mathcal{B}_{\lambda}^{g}}P\left(  x,\xi\right)  \left(
g\left(  \xi_{\lambda}\right)  ^{q/r}-g\left(  \xi\right)  ^{q/r}\right)
^{r}d\xi\right)  ^{1/r}\\
& \leq\frac{2q}{r}\left(  \int_{\mathcal{B}_{\lambda}^{g}}P\left(
x,\xi\right)  g\left(  \xi_{\lambda}\right)  ^{q-r}\left(  g\left(
\xi_{\lambda}\right)  -g\left(  \xi\right)  \right)  ^{r}d\xi\right)  ^{1/r}.
\end{align*}
Hence%
\begin{align*}
& \left\vert v_{\lambda}^{1/r}-v^{1/r}\right\vert _{L^{\frac{n\left(
q+1\right)  r}{\left(  n-1\right)  q}}\left(  \mathcal{B}_{\lambda}%
^{v}\right)  }\\
& \leq\frac{2q}{r}\left\vert \int_{\mathcal{B}_{\lambda}^{g}}P\left(
x,\xi\right)  g\left(  \xi_{\lambda}\right)  ^{q-r}\left(  g\left(
\xi_{\lambda}\right)  -g\left(  \xi\right)  \right)  ^{r}d\xi\right\vert
_{L^{\frac{n\left(  q+1\right)  }{\left(  n-1\right)  q}}\left(
\mathcal{B}_{\lambda}^{v}\right)  }^{1/r}\\
& \leq c\left(  n,q,r\right)  \left\vert g_{\lambda}^{q-r}\left(  g_{\lambda
}-g\right)  ^{r}\right\vert _{L^{\frac{q+1}{q}}\left(  \mathcal{B}_{\lambda
}^{g}\right)  }^{1/r}\\
& \leq c\left(  n,q,r\right)  \left\vert g_{\lambda}\right\vert _{L^{q+1}%
\left(  \mathcal{B}_{\lambda}^{g}\right)  }^{\frac{q-r}{r}}\left\vert
g_{\lambda}-g\right\vert _{L^{q+1}\left(  \mathcal{B}_{\lambda}^{g}\right)  }.
\end{align*}
Combining the two inequalities together we see%
\[
\left\vert g_{\lambda}-g\right\vert _{L^{q+1}\left(  \mathcal{B}_{\lambda}%
^{g}\right)  }\leq c\left(  n,q,r\right)  \left\vert v\right\vert
_{L^{\frac{n\left(  q+1\right)  }{\left(  n-1\right)  q}}\left(
\mathbb{R}_{+}^{n}\right)  }^{\frac{q+n}{\left(  n-1\right)  q}-r^{-1}%
}\left\vert g\right\vert _{L^{q+1}\left(  2\lambda e_{1}^{\prime}%
-\mathcal{B}_{\lambda}^{g}\right)  }^{\frac{q-r}{r}}\left\vert g_{\lambda
}-g\right\vert _{L^{q+1}\left(  \mathcal{B}_{\lambda}^{g}\right)  }.
\]
With this inequality at hand, we may proceed in the same way as in the Case
\ref{case6.1} to get the conclusion that $g$ is radial symmetric with respect
to some point and strictly decreasing along the radial direction.
\end{proof}

Next we look at the special power $p=\frac{2\left(  n-1\right)  }{n-2}$.

\begin{proof}
[Proof of Proposition \ref{prop6.1}]If we know $f\in L^{\frac{2\left(
n-1\right)  }{n-2}}\left(  \mathbb{R}^{n-1}\right)  $, then it follows from
Theorem \ref{thm6.1} that $f\in C^{\infty}\left(  \mathbb{R}^{n-1}\right)  $,
it is strictly positive and radial symmetric with respect to some point. By
translation we may assume $f$ is radial symmetric with respect to $0$.

On the other hand, if $f$ is a solution to the equation, let $u\left(
x\right)  =\left(  Pf\right)  \left(  x\right)  $, $\widetilde{f}\left(
\xi\right)  =\frac{1}{\left\vert \xi\right\vert ^{n-2}}f\left(  \frac{\xi
}{\left\vert \xi\right\vert ^{2}}\right)  $ and $\widetilde{u}\left(
x\right)  =\frac{1}{\left\vert x\right\vert ^{n-2}}u\left(  \frac
{x}{\left\vert x\right\vert ^{2}}\right)  $, by change of variable we know%
\[
\widetilde{u}\left(  x\right)  =\left(  P\widetilde{f}\right)  \left(
x\right)  ,\quad\widetilde{f}\left(  \xi\right)  ^{\frac{n}{n-2}}%
=\int_{\mathbb{R}_{+}^{n}}P\left(  x,\xi\right)  \widetilde{u}\left(
x\right)  ^{\frac{n+2}{n-2}}dx
\]
and $\left\vert \widetilde{f}\right\vert _{L^{\frac{2\left(  n-1\right)
}{n-2}}\left(  \mathbb{R}^{n-1}\right)  }=\left\vert f\right\vert
_{L^{\frac{2\left(  n-1\right)  }{n-2}}\left(  \mathbb{R}^{n-1}\right)  }$. In
particular, $\widetilde{f}\in L^{\frac{2\left(  n-1\right)  }{n-2}}$ and
satisfies the same equation.

Let $e_{1}=\left(  1,0,\cdots,0\right)  \in\mathbb{R}^{n}$, then it follows
from Theorem \ref{thm6.1} that $f_{1}\left(  \xi\right)  =\frac{1}{\left\vert
\xi\right\vert ^{n-2}}f\left(  \frac{\xi}{\left\vert \xi\right\vert ^{2}%
}-e_{1}^{\prime}\right)  $ is smooth and radial symmetric with respect to some
point. It follows from Proposition \ref{prop4.1} and the fact that $f\in
L^{\frac{2\left(  n-1\right)  }{n-2}}\left(  \mathbb{R}^{n-1}\right)  $ that
for some $c_{1}>0$ and $c_{2}>0$, $f\left(  \xi\right)  =\left(
c_{1}\left\vert \xi\right\vert ^{2}+c_{2}\right)  ^{-\frac{n-2}{2}}$. Since
$f$ satisfies the equation, it follows that for some $\lambda>0$, $f\left(
\xi\right)  =c\left(  n\right)  \left(  \frac{\lambda}{\lambda^{2}+\left\vert
\xi\right\vert ^{2}}\right)  ^{\frac{n-2}{2}}$.

Next we want to show under the assumption of the Proposition \ref{prop6.1},
$f$ always lies in $L^{\frac{2\left(  n-1\right)  }{n-2}}\left(
\mathbb{R}^{n-1}\right)  $. This will be proved by contradiction. Indeed, if
this is not the case, then $\int_{\mathbb{R}^{n-1}}f\left(  \xi\right)
^{\frac{2\left(  n-1\right)  }{n-2}}d\xi=\infty$. Let $g_{0}\left(
\xi\right)  =f\left(  \xi\right)  ^{\frac{n}{n-2}}$, $v_{0}\left(  x\right)
=\left(  Pf\right)  \left(  x\right)  $, then $g_{0}\in L_{loc}^{\frac
{2\left(  n-1\right)  }{n}}\left(  \mathbb{R}^{n-1}\right)  $, $\int
_{\mathbb{R}^{n-1}}g_{0}\left(  \xi\right)  ^{\frac{2\left(  n-1\right)  }{n}%
}d\xi=\infty$ and%
\[
v_{0}\left(  x\right)  =\int_{\mathbb{R}^{n-1}}P\left(  x,\xi\right)
g_{0}\left(  \xi\right)  ^{\frac{n-2}{n}}d\xi,\quad g_{0}\left(  \xi\right)
=\int_{\mathbb{R}_{+}^{n}}P\left(  x,\xi\right)  v_{0}\left(  x\right)
^{\frac{n+2}{n-2}}dx.
\]
It follows from the proof of Theorem \ref{thm5.1} that $g_{0}\in C^{\infty
}\left(  \mathbb{R}^{n-1}\right)  $ and $v_{0}\in C^{\infty}\left(
\overline{\mathbb{R}}_{+}^{n}\right)  $.

Let $g\left(  \xi\right)  =\frac{1}{\left\vert \xi\right\vert ^{n}}%
g_{0}\left(  \frac{\xi}{\left\vert \xi\right\vert ^{2}}\right)  $, $v\left(
x\right)  =\frac{1}{\left\vert x\right\vert ^{n-2}}v_{0}\left(  \frac
{x}{\left\vert x\right\vert ^{2}}\right)  $, then%
\[
v\left(  x\right)  =\int_{\mathbb{R}^{n-1}}P\left(  x,\xi\right)  g\left(
\xi\right)  ^{\frac{n-2}{n}}d\xi,\quad g\left(  \xi\right)  =\int
_{\mathbb{R}_{+}^{n}}P\left(  x,\xi\right)  v\left(  x\right)  ^{\frac
{n+2}{n-2}}dx.
\]
Moreover for any $R>0$, $\int_{\mathbb{R}^{n-1}\backslash B_{R}}g\left(
\xi\right)  ^{\frac{2\left(  n-1\right)  }{n}}d\xi<\infty$ and $\int
_{\mathbb{R}^{n-1}}g\left(  \xi\right)  ^{\frac{2\left(  n-1\right)  }{n}}%
d\xi=\infty$. For $\lambda>0$, we define $H_{\lambda},g_{\lambda}$ as in the
Case \ref{case6.1} of the proof of Theorem \ref{thm6.1}, but let
$\mathcal{B}_{\lambda}^{g}=\left\{  \xi\in H_{\lambda}\backslash\left\{
0\right\}  :g_{\lambda}\left(  \xi\right)  >g\left(  \xi\right)  \right\}  $.
Put the number in the proof of Theorem \ref{thm6.1} $r=\frac{n+2}{n-2}$, then
the same argument shows%
\[
\left\vert g_{\lambda}^{\frac{n-2}{n+2}}-g^{\frac{n-2}{n+2}}\right\vert
_{L^{\frac{2\left(  n-1\right)  \left(  n+2\right)  }{n\left(  n-2\right)  }%
}\left(  \mathcal{B}_{\lambda}^{g}\right)  }\leq c\left(  n\right)  \left\vert
g\right\vert _{L^{\frac{2\left(  n-1\right)  }{n}}\left(  2\lambda
e_{1}^{\prime}-\mathcal{B}_{\lambda}^{g}\right)  }^{\frac{2\left(  n-2\right)
}{n\left(  n+2\right)  }}\left\vert g_{\lambda}^{\frac{n-2}{n+2}}%
-g^{\frac{n-2}{n+2}}\right\vert _{L^{\frac{2\left(  n-1\right)  \left(
n+2\right)  }{n\left(  n-2\right)  }}\left(  \mathcal{B}_{\lambda}^{g}\right)
}.
\]
Note that for $\xi\in\mathcal{B}_{\lambda}^{g}$, $g_{\lambda}\left(
\xi\right)  >g\left(  \xi\right)  $, hence%
\[
\int_{\mathcal{B}_{\lambda}^{g}}g\left(  \xi\right)  ^{\frac{2\left(
n-1\right)  }{n}}d\xi\leq\int_{\mathcal{B}_{\lambda}^{g}}g_{\lambda}\left(
\xi\right)  ^{\frac{2\left(  n-1\right)  }{n}}d\xi\leq\int_{\mathbb{R}%
^{n-1}\backslash H_{\lambda}}g\left(  \xi\right)  ^{\frac{2\left(  n-1\right)
}{n}}d\xi<\infty.
\]
When $\lambda$ is large enough, it implies%
\[
\left\vert g_{\lambda}^{\frac{n-2}{n+2}}-g^{\frac{n-2}{n+2}}\right\vert
_{L^{\frac{2\left(  n-1\right)  \left(  n+2\right)  }{n\left(  n-2\right)  }%
}\left(  \mathcal{B}_{\lambda}^{g}\right)  }\leq\frac{1}{2}\left\vert
g_{\lambda}^{\frac{n-2}{n+2}}-g^{\frac{n-2}{n+2}}\right\vert _{L^{\frac
{2\left(  n-1\right)  \left(  n+2\right)  }{n\left(  n-2\right)  }}\left(
\mathcal{B}_{\lambda}^{g}\right)  }%
\]
and hence $\left\vert g_{\lambda}^{\frac{n-2}{n+2}}-g^{\frac{n-2}{n+2}%
}\right\vert _{L^{\frac{2\left(  n-1\right)  \left(  n+2\right)  }{n\left(
n-2\right)  }}\left(  \mathcal{B}_{\lambda}^{g}\right)  }=0$, $\mathcal{B}%
_{\lambda}^{g}=\emptyset$. Let
\[
\lambda_{0}=\inf\left\{  \lambda>0:\mathcal{B}_{\lambda^{\prime}}%
^{g}=\emptyset\text{ for all }\lambda^{\prime}\geq\lambda\right\}  .
\]
We claim $\lambda_{0}=0$. Indeed if this is not the case, then $\lambda_{0}%
>0$. We may argue as in the Case \ref{case6.1} of the proof of Theorem
\ref{thm6.1} and get $g_{\lambda_{0}}=g$. In particular, this would imply
$\int_{\mathbb{R}^{n-1}}g\left(  \xi\right)  ^{\frac{2\left(  n-1\right)  }%
{n}}d\xi<\infty$, a contradiction. It follows that $\lambda_{0}=0$ and
$g\left(  \xi_{1},\xi^{\prime\prime}\right)  \geq g\left(  -\xi_{1}%
,\xi^{\prime\prime}\right)  $ for $\xi_{1}<0$. Since we may perform this
process along any direction, we see $g$ must be radial symmetric with respect
to $0$. Hence $g_{0}$ must be radial symmetric with respect to $0$. For any
$\zeta\in\mathbb{R}^{n-1}$, we may apply the argument to $g_{0}\left(
\cdot+\zeta\right)  $ and deduce that $g_{0}$ is also radial symmetric with
respect to $\zeta$, hence $g_{0}$ must be a constant function, so is $f$. But
this contradicts with the fact that $f$ satisfies the equation.
\end{proof}


\begin{thebibliography}{9999}                                                                                             %
\bibitem[C]{C}T. Carleman. Zur Theorie de Minimalfl\"{a}chen. \textit{Math Z.
}\textbf{9} (1921), 154--160.

\bibitem[CGS]{CGS}L. A. Caffarelli, B. Gidas and J. Spruck. Asymptotic
symmetry and local behavior of semilinear elliptic equations with critical
Sobolev growth. \textit{Comm Pure Appl Math.} \textbf{42} (1989), no. 3, 271--297.

\bibitem[CL]{CL}E. A. Carlen and M. Loss. Extremals of functionals with
competing symmetries. \textit{J Funct Anal}. \textbf{88} (1990), no. 2, 437--456.

\bibitem[ChL]{ChL}W. X. Chen and C. M. Li. Regularity of solutions for a
system of integral equations. \textit{Comm Pure Appl Anal}, \textbf{4}(1),
(2005), 1--8.

\bibitem[CLO1]{CLO1}W. X. Chen, C. M. Li and B. Ou. Classification of
solutions for a system of integral equations. \textit{Comm in Partial
Differential Equations}, \textbf{30} (2005), 59--65.

\bibitem[CLO2]{CLO2}W. X. Chen, C. M. Li and B. Ou. Classification of
solutions for an integral equation. \textit{Comm Pure Appl Math}, \textbf{59}
(2006), no. 3, 330--343.

\bibitem[CLO3]{CLO3}W. X. Chen, C. M. Li and B. Ou. Alternative proofs on the
radial symmetry and monotonicity of positive regular solutions to a singular
integral equation. IMA Preprint Series \#2044.

\bibitem[EG]{EG}L. C. Evans and R. F. Gariepy. Measure theory and fine
properties of functions. Studies in Advanced Mathematics. CRC Press, Boca
Raton, FL, 1992.

\bibitem[GNN]{GNN}B. Gidas, W. M. Ni and L. Nirenberg. Symmetry and related
properties via the maximum principle. \textit{Communications in Mathematical
Physics}. \textbf{68} (1979), 209--243.

\bibitem[GT]{GT}D. Gilbarg \ and N. S. Trudinger. Elliptic partial
differential equations of second order. 2nd edition, 3rd printing. Berlin:
Springer-Verlag, 1998.

\bibitem[HL]{HL}Q. Han and F. H. Lin. Elliptic partial differential equations.
Courant Lecture Notes, volume 1. American Mathematical Society, 2000.

\bibitem[Hn]{Hn}F. B. Hang. On the integral systems related to
Hardy-Littlewood-Sobolev inequality. Preprint, 2005.

\bibitem[HWY]{HWY}F. B. Hang, X. D. Wang and X. D. Yan. An integral equation
in conformal geometry. In preparation.

\bibitem[L]{L}Y. Y. Li. Remark on some conformally invariant integral
equations: the method of moving spheres. \textit{Journal of European
Mathematical Society}. \textbf{6} (2004), 153--180.

\bibitem[LZ]{LZ}Y. Li and M. Zhu. Uniqueness theorems through the method of
moving spheres. \textit{Duke Mathematical Journal}.\textit{\ }\textbf{80 }no.
2 (1995), 383--417.

\bibitem[Li1]{Li1}E. H. Lieb. Existence and uniqueness of the minimizing
solution of Choquard's nonlinear equation. \textit{Studies in Appl Math.}
\textbf{57} (1977), 93--105.

\bibitem[Li2]{Li2}E. H. Lieb. Sharp constants in the Hardy-Littlewood-Sobolev
and related inequalities. \textit{Ann of Math}, \textbf{118}(2), (1983), 349--374.

\bibitem[LiL]{LiL}E. H. Lieb and M. Loss. Analysis. Second edition. Graduate
Studies in Mathematics, 14. American Mathematical Society, Providence, RI, 2001.

\bibitem[Lion]{Lion}P. L. Lions. The concentration-compactness principle in
the calculus of variations. The limit case. II. \textit{Rev. Mat.
Iberoamericana} \textbf{1} (1985), no. 2, 45--121.

\bibitem[O]{O}B. Ou. A remark on a singular integral equation. \textit{Houston
J Math}. \textbf{25} (1999), 181--184.

\bibitem[S]{S}E. M. Stein. Singular integrals and differentiability properties
of functions. Princeton University Press, Princeton, New Jersey, 1970.

\bibitem[SW]{SW}E. M. Stein and G. Weiss. Introduction to Fourier analysis on
Euclidean spaces. Princeton Mathematical Series, No. 32. Princeton University
Press, Princeton NJ, 1971.
\end{thebibliography}
\end{document}